\documentclass[12pt]{amsart}
\usepackage{latexsym,amssymb,amsrefs,graphicx}
\usepackage[all]{xy}

\newtheorem{lemma}{Lemma}

\newtheorem{theorem}[lemma]{Theorem}
\newtheorem{corollary}[lemma]{Corollary}
\newtheorem{definition}[lemma]{Definition}

\theoremstyle{definition}

\newcommand{\proofstart}{\smallskip\noindent Proof:\ }

\newcommand{\proofend}{\begin{flushright}
                       $\Box$
                       \end{flushright}}

\begin{document}

\title[Cascades and perturbed Morse-Bott functions]{Cascades and perturbed Morse-Bott functions}

\author{Augustin Banyaga}
\address{Department of Mathematics \\
         Penn State University \\
         University Park, PA 16802}
\email{banyaga@math.psu.edu}

\author{David E. Hurtubise}
\address{Department of Mathematics and Statistics\\
         Penn State Altoona\\
         Altoona, PA 16601-3760}
\email{Hurtubise@psu.edu}

\subjclass[2000]{Primary: 57R70 Secondary: 37D05 37D15 58E05}

\begin{abstract}
Let $f:M \rightarrow \mathbb{R}$ be a Morse-Bott function on a finite dimensional closed smooth
manifold $M$. Choosing an appropriate Riemannian metric on $M$ and Morse-Smale
functions $f_j:C_j \rightarrow \mathbb{R}$ on the critical submanifolds $C_j$, one can
construct a Morse chain complex whose boundary operator is defined by
counting cascades \cite{FraTheA}. Similar data, which also includes a parameter 
$\varepsilon > 0$ that scales the Morse-Smale functions $f_j$, can be used to define an
explicit perturbation of the Morse-Bott function $f$ to a Morse-Smale function 
$h_\varepsilon:M \rightarrow \mathbb{R}$ \cite{AusMor} \cite{BanDyn}. In this paper we
show that the Morse-Smale-Witten chain complex of $h_\varepsilon$ is the same as the Morse 
chain complex defined using cascades for any $\varepsilon >0$ sufficiently small. That is,
the two chain complexes have the same generators, and their boundary operators are the same
(up to a choice of sign). Thus, the Morse Homology Theorem implies that the homology of
the cascade chain complex of $f:M \rightarrow \mathbb{R}$ is isomorphic to the singular
homology $H_\ast(M;\mathbb{Z})$.
\end{abstract}

\maketitle


\section{Introduction}

Let $f:M\rightarrow \mathbb{R}$ be a Morse-Bott function on a finite dimensional closed
smooth Riemannian manifold $(M,g)$ with connected critical submanifolds $C_j$ for
$j=1, \ldots , l$. There are at least three approaches to computing the homology
of $M$ using moduli spaces of gradient flow lines:
\begin{enumerate}
  \item Perturb $f:M \rightarrow \mathbb{R}$ to a Morse-Smale function and use the
        Morse-Smale-Witten chain complex, whose boundary operator is defined using
        moduli spaces of gradient flow lines of the perturbed function (see for instance
        \cite{BanLec}, \cite{SchMor}, and the references therein).
  \item Introduce Morse functions $f_j:C_j \rightarrow \mathbb{R}$ on the critical
        submanifolds $C_1,\ldots, C_l$ and use a Morse chain complex whose boundary 
        operator is defined using moduli spaces of cascades \cite{FraTheA}. 
  \item Use the Morse-Bott-Smale multicomplex, where the homomorphisms in the multicomplex
        are defined using fibered products of moduli spaces of gradient flow lines
        of the Morse-Bott function $f:M\rightarrow \mathbb{R}$ \cite{BanMor}.
\end{enumerate}
A fourth approach might involve using the filtration determined by the Morse-Bott
function $f:M \rightarrow \mathbb{R}$ to define a spectral sequence, but the
differentials in the spectral sequence determined by the filtration are not defined
using moduli spaces of gradient flow lines (see \cite{BanMor} and \cite{HurMul}). In
addition, there are approaches to computing the cohomology/homology of $M$ from a Morse-Bott 
function using differential forms and/or currents \cite{AusMor} \cite{ChoMor} \cite{LatGra},
but we will not discuss differential forms or currents in this paper.

\vspace{.1 in}

The main goal of this paper is to show that for a finite dimensional closed smooth manifold
$M$ the first two approaches are essentially the same. That is, the auxiliary Morse functions
$f_j:C_j \rightarrow \mathbb{R}$ on the critical submanifolds $C_j$ for $j=1, \ldots ,l$ 
required to define the cascade chain complex and a parameter $\varepsilon > 0$ determine 
an explicit perturbation of the Morse-Bott function $f:M \rightarrow \mathbb{R}$ to a Morse
function $h_\varepsilon:M \rightarrow \mathbb{R}$ \cite{AusMor} \cite{BanDyn}. Moreover, 
under certain transversality assumptions the Morse-Smale-Witten chain complex of 
$h_\varepsilon:M \rightarrow \mathbb{R}$ has the same generators and the same boundary
operator as the cascade chain complex (up to a choice of sign).

\vspace{.1 in}

We now describe the cascade chain complex for a Morse-Bott function. To the best of our
knowledge, moduli spaces of cascades were first introduced within the context of symplectic
Floer homology \cite{FraTheA}. However, moduli spaces of cascades have since been used in the
contexts of contact homology and gauge theory \cite{BouAne} \cite{BouSym} \cite{CieAFl}
\cite{SwoMor}. Our approach to constructing moduli spaces of cascades and their
compactifications is given in Sections \ref{cascadesection} and \ref{broken} for a function 
$f:M\rightarrow \mathbb{R}$ on a finite dimensional closed smooth Riemannian manifold 
$(M,g)$ that satisfies the Morse-Bott-Smale transversality condition. The moduli spaces of
cascades are constructed using finite dimensional fibered products similar to those found
in \cite{BanMor}, and the compactifications of the moduli spaces are described in terms 
of the Hausdorff topology.


\subsection*{Cascades}

Let $f:M \rightarrow \mathbb{R}$ be a Morse-Bott function on a finite dimensional closed
smooth Riemannian manifold $(M,g)$ with connected critical submanifolds $C_1, \ldots, C_l$.
Choose Morse-Smale functions $f_j:C_j \rightarrow \mathbb{R}$ on the critical submanifolds
for all $j=1,\ldots ,l$, and define the \textit{total index} of a critical point of $f_j$ to
be its Morse index on $C_j$ plus the Morse-Bott index of the critical submanifold
$C_j$. Roughly speaking, a \textit{cascade} between two critical points is a concatenation
of some gradient flow lines of the function $f$ and pieces of the gradient flow lines of the
functions $f_j$ on the critical submanifolds. Choosing appropriate Riemannian metrics on $M$
and the critical submanifolds $C_j$ it is shown in the appendix to \cite{FraTheA} that the
moduli space of cascades $\mathcal{M}^c(q,p)$ between two critical points $q$ and $p$ is a
smooth manifold of dimension $\lambda_q - \lambda_p - 1$, where $\lambda_q$ and $\lambda_p$
denote the total indices of $q$ and $p$ respectively. Moreover, $\mathcal{M}^c(q,p)$ has a
compactification consisting of broken flow lines with cascades between $q$ and $p$.

Since the moduli space of cascades $\mathcal{M}^c(q,p)$ has properties similar to those of
a moduli space of gradient flow lines of a Morse-Smale function, it is natural to define a chain
complex analogous to the Morse-Smale-Witten chain complex but using moduli spaces of cascades
in place of moduli spaces of gradient flow lines. Thus, we define the $k^{\text{th}}$ chain
group $C_k^c(f)$ to be the free abelian group generated by the critical points of total index
$k$ of the Morse-Smale functions $f_j$ for all $j=1,\ldots ,l$. In the appendix to \cite{FraTheA}
a boundary operator $\partial^c$ is defined by counting the number of cascades between
critical points of relative index one mod $2$, and a continuation theorem is stated that implies
that the homology of the chain complex $(C_\ast^c(f)\otimes \mathbb{Z}_2, \partial^c)$ is 
isomorphic to the singular homology $H_\ast(M;\mathbb{Z}_2)$. In Section \ref{correspondencesection}
of this paper we show that it is possible to define the boundary operator $\partial^c$ over 
$\mathbb{Z}$ by counting the elements of $\mathcal{M}^c(q,p)$ with sign when $\lambda_q - \lambda_p = 1$,
and we prove that the homology of the resulting chain complex $(C_\ast^c(f), \partial^c)$ is 
isomorphic to the singular homology $H_\ast(M;\mathbb{Z})$. 


\subsection*{Perturbing the Morse-Bott function}

The Morse-Smale functions $f_j:C_j \rightarrow \mathbb{R}$ chosen to define the
chain complex $(C_\ast^c(f),\partial^c)$ can also be used to define an explicit perturbation
of the Morse-Bott function $f:M\rightarrow \mathbb{R}$ to a Morse-Smale function 
$h_\varepsilon:M\rightarrow \mathbb{R}$. This perturbation technique was used in \cite{AusMor}
in relation to a de Rham version of Morse-Bott cohomology. It was also used in \cite{BanDyn}
to give a dynamical systems approach to the proof of the Morse-Bott inequalities with
somewhat different orientation assumptions than the classical ``half-space'' method using
the Thom Isomorphism Theorem (see \cite{BotMor}, Appendix C of \cite{FarTop}, and Section 2.6 
of \cite{NicAnI}).

\vspace{.1 in}

To define the Morse-Smale function $h_\varepsilon:M\rightarrow \mathbb{R}$ near $f$ choose 
``small'' tubular neighborhoods $T_j$ of each of the critical submanifolds $C_j$ for all 
$j=1,\ldots ,l$ and extend the Morse-Smale functions $f_j$ to the tubular neighborhoods $T_j$
by making them constant in the direction normal to $C_j$. Choose bump functions $\rho_j$ on the
tubular neighborhoods $T_j$ for all $j=1,\ldots ,l$ that are equal to one in an open neighborhood
of $C_j$, constant in the direction parallel to $C_j$, and equal to zero outside of $T_j$. 
The function
$$
h_\varepsilon = f + \varepsilon \left( \sum_{j=1}^l \rho_j f_j \right)
$$
is a Morse function near $f$ for any sufficiently small $\varepsilon >0$, and the critical
set of $h_\varepsilon$ is the union of the critical points of the functions
$f_j:C_j \rightarrow \mathbb{R}$ for $j=1, \ldots, l$.  In fact, the total index $\lambda_q$
of a critical point $q$ is the same as the Morse index of $q$ viewed as a critical point of
$h_\varepsilon:M \rightarrow \mathbb{R}$.


\subsection*{Correspondence}

If we choose the Riemannian metric $g$ on $M$ so that $h_\varepsilon:M \rightarrow
\mathbb{R}$ satisfies the Morse-Smale transversality condition with respect to $g$, then
the moduli space of gradient flow lines of $h_\varepsilon$ between two critical points $q$
and $p$ is a smooth manifold with $\text{dim }\mathcal{M}_{h_\varepsilon}(q,p) = 
\lambda_q - \lambda_p - 1$. We show in Section \ref{cascadesection} that if $f:M \rightarrow 
\mathbb{R}$ satisfies the Morse-Bott-Smale transversality condition and
we choose the Morse functions $f_j$ on the critical submanifolds so that some additional
transversality conditions are satisfied, then the moduli space of cascades $\mathcal{M}^c(q,p)$
is also a smooth manifold of dimension $\lambda_q - \lambda_p - 1$. 

\vspace{.1 in}
In Section \ref{correspondencesection} we prove that when the dimension of these moduli spaces
is zero they have the same number of elements.

\begin{theorem}[Correspondence of Moduli Spaces]
Let $p,q \in \text{Cr}(h_\varepsilon)$ with $\lambda_q - \lambda_p = 1$.
For any sufficiently small $\varepsilon>0$ there is a bijection between unparameterized
cascades and unparameterized gradient flow lines of the Morse-Smale function
$h_\varepsilon:M \rightarrow \mathbb{R}$ between $q$ and $p$,
$$
\mathcal{M}^c(q,p) \leftrightarrow \mathcal{M}_{h_\varepsilon}(q,p).
$$
\end{theorem}

\noindent
Choosing orientations on the unstable manifolds of the Morse-Smale function
$h_\varepsilon:M \rightarrow \mathbb{R}$ associates a sign $\pm 1$ to each component of 
$\mathcal{M}_{h_\varepsilon}(q,p)$ when $\lambda_q - \lambda_p = 1$,
and thus we can use the correspondence theorem for moduli spaces to transport the signs
to the components of $\mathcal{M}^c(q,p)$. This allows us to define the boundary operator
in the cascade chain complex over $\mathbb{Z}$, and we have the following as an immediate
corollary.

\begin{corollary}[Correspondence of Chain Complexes]
For $\varepsilon >0$ sufficiently small, the Morse-Smale-Witten chain complex
$(C_\ast(h_\varepsilon),\partial)$ associated to the perturbation
$$
h_\varepsilon = f + \varepsilon \left( \sum_{j=1}^l \rho_j f_j \right)
$$
of a Morse-Bott function $f:M \rightarrow \mathbb{R}$ is the same as the cascade
chain complex $(C_\ast^c(f),\partial^c)$. That is, the chain groups of both
complexes have the same generators and their boundary operators are the same 
(up to a choice of sign).
\end{corollary}

\smallskip\noindent
This corollary, together with the Morse Homology Theorem, implies immediately that the
homology of the chain complex $(C_\ast^c(f),\partial^c)$ is isomorphic to the singular
homology $H_\ast(M;\mathbb{Z})$.


\subsection*{Outline of the paper}

In Section \ref{MorseChain} we recall some basic definitions and facts about
the Morse-Smale-Witten chain complex.  In Section \ref{cascadesection} we give
a detailed construction of the smooth moduli space of cascades $\mathcal{M}^c(q,p)$
under the assumption that $f:M \rightarrow \mathbb{R}$ satisfies the
Morse-Bott-Smale transversality condition with respect to the metric $g$ on $M$.
Our construction requires that the Morse functions $f_j:C_j \rightarrow \mathbb{R}$ satisfy 
the Morse-Smale transversality condition with respect to the restriction of the
Riemannian metric $g$ to the critical submanifolds for all $j=1, \ldots, l$ and
that all the unstable and stable manifolds on the critical submanifolds are transverse
to certain beginning and endpoint maps (Definition \ref{beginningendtransverse}).
Lemma \ref{perturb} shows that it is always possible to choose the auxiliary Morse
functions $f_j:M \rightarrow \mathbb{R}$ so that these transversality conditions are
satisfied.  Theorem \ref{manifold} shows that under the above assumptions
$\mathcal{M}^c(q,p)$ is a smooth manifold of dimension $\lambda_q - \lambda_p - 1$
that is stratified by smooth manifolds with corners.

In Section \ref{broken} we study the compactness properties of $\mathcal{M}^c(q,p)$.
We show using the Hausdorff metric that $\mathcal{M}^c(q,p)$ can be compactified
using broken flow lines with cascades, which implies that $\mathcal{M}^c(q,p)$ is
compact when $\lambda_q - \lambda_p = 1$. In Section \ref{correspondencesection}
we give a detailed construction of the perturbation $h_\varepsilon:M \rightarrow \mathbb{R}$,
and we prove that it is possible to choose a single Riemannian metric $g$ so that
$h_\varepsilon:M \rightarrow \mathbb{R}$ satisfies the Morse-Smale transversality
condition with respect to $g$ for all $\varepsilon > 0$ sufficiently small 
(Lemma \ref{perturbmetric}).  We also prove that as $\varepsilon \rightarrow 0$
a sequence of gradient flow lines of $h_\varepsilon$ between two critical points
$q$ and $p$ must have a subsequence that converges to a broken flow line with cascades 
from $q$ to $p$ (Lemma \ref{degenerating}).

The correspondence theorem for moduli spaces (Theorem \ref{correspondence}) is proved
in Section \ref{correspondencesection} using recent results from geometric singular
perturbation theory.  In particular, our proof uses the Exchange Lemma for fast-slow systems
\cite{JonGen} \cite{SchExcI} \cite{SchExcII} which says (roughly) that a manifold $M_0$ that is
transverse to the stable manifold of a normally hyperbolic locally invariant submanifold
$C$ will have subsets that flow forward in time under the full fast-slow system to be
near subsets of the unstable manifold of $C$.  The correspondence theorem for the
Morse-Smale-Witten chain complex of $h_\varepsilon:M \rightarrow \mathbb{R}$ and the
cascade chain complex (Corollary \ref{chaincorrespondence}) follows as an
immediate corollary to the correspondence theorem for moduli spaces.


\section{The Morse-Smale-Witten chain complex}\label{MorseChain}

In this section we briefly recall the construction of the Morse-Smale-Witten
chain complex and the Morse Homology Theorem.  For more details see \cite{BanLec}.

\smallskip
Let $Cr(f) = \{p \in M |\, df_p = 0 \}$ denote the set of critical points of a 
smooth function $f:M \rightarrow \mathbb{R}$ on a  smooth $m$-dimensional manifold
$M$. A critical point $p \in Cr(f)$ is said to be \textbf{nondegenerate} if
and only if the Hessian $H_p(f)$ is nondegenerate.  The \textbf{index} $\lambda_p$
of a nondegenerate critical point $p$ is the dimension of the subspace of $T_pM$
where $H_p(f)$ is negative definite. If all the critical points of $f$ are
non-degenerate, then $f$ is called a \textbf{Morse function}.

If $f:M \rightarrow \mathbb{R}$ is a Morse function on a finite dimensional 
compact smooth Riemannian manifold $(M,g)$, then the \textbf{stable manifold} $W^s_f(p)$
and the \textbf{unstable manifold} $W^u_f(p)$ of a critical point $p \in Cr(f)$ are
defined to be
\begin{eqnarray*}
W^s_f(p) & = & \{ x\in M | \lim_{t \rightarrow \infty} \varphi_t(x) = p \}\\
W^u_f(p) & = & \{ x\in M | \lim_{t \rightarrow -\infty} \varphi_t(x) = p \}
\end{eqnarray*}
where $\varphi_t$ is the 1-parameter group of diffeomorphisms generated by
minus the gradient vector field, i.e. $-\nabla f$. The index of $p$ coincides 
with the dimension of $W^u_f(p)$. The Stable/Unstable Manifold Theorem for a Morse Function
says that the tangent space at $p$ splits as 
$$
T_pM = T^s_pM \oplus T_p^uM
$$
where the Hessian is positive definite on $T_p^sM \stackrel{def}{=} T_p W^s_f(p)$ and 
negative definite on $T_p^uM \stackrel{def}{=} T_p W^u_f(p)$.  Moreover, the
stable and unstable manifolds of $p$ are surjective images of smooth 
embeddings
\begin{eqnarray*}
E^s: T_p^sM & \rightarrow & W^s_f(p) \subseteq M\\
E^u: T_p^uM & \rightarrow & W^u_f(p) \subseteq M. 
\end{eqnarray*}
Hence, $W^s_f(p)$ is a smoothly embedded open disk of dimension $m - \lambda_p$, 
and $W^u_f(p)$ is a smoothly embedded open disk of dimension $\lambda_p$.

If the stable and unstable manifolds of a Morse function $f:M \rightarrow
\mathbb{R}$ all intersect transversally, then the function $f$ is called
\textbf{Morse-Smale}. For any metric $g$ on $M$ the set of smooth Morse-Smale
functions is dense by the Kupka-Smale Theorem (Theorem 6.6 and Remark 6.7
of \cite{BanLec}), and for a given Morse function $f:M \rightarrow \mathbb{R}$ one can
choose a Riemannian metric on $M$ so that $f$ is Morse-Smale with respect to the chosen
metric (Theorem 2.20 of \cite{AbbLec}). Moreover, if $f$ is Morse-Smale and $p,q\in Cr(f)$
then $W_f(q,p) = W^u_f(q) \cap W^s_f(p)$ is an embedded submanifold of $M$ of dimension
$\lambda_q  - \lambda_p$, and when $\lambda_q - \lambda_p = 1$ the number of gradient flow
lines from $q$ to $p$ is finite (Corollary 6.29 of \cite{BanLec}).

If we choose an orientation for each of the unstable manifolds of $f$, then there is an
induced orientation on the normal bundles of the stable manifolds. Thus, we can define an integer
associated to any two critical points $p$ and $q$ of relative index one by counting the number
of gradient flow lines from $q$ to $p$ with signs determined by the orientations. This integer
is denoted by $n_f(q,p) = \# \mathcal{M}_f(q,p)$, where $\mathcal{M}_f(q,p) = 
W_f(q,p)/\mathbb{R}$ is the moduli space of gradient flow lines of $f$ from $q$ to $p$. 
The \textbf{Morse-Smale-Witten chain complex} is defined to be the chain complex 
$(C_\ast(f),\partial_\ast)$ where $C_k(f)$ is the free abelian group generated by the 
critical points $q$ of index $k$ and the boundary operator 
$\partial_k:C_k(f) \rightarrow C_{k-1}(f)$ is given by
$$
\partial_k(q)\ \ = \sum_{p \in Cr_{k-1}(f)} n_f(q,p)p.
$$

\begin{theorem}[Morse Homology Theorem]\label{Morsehomology}
The pair $(C_\ast(f),\partial_\ast)$ is a chain complex, and the homology
of $(C_\ast(f),\partial_\ast)$ is isomorphic to the singular homology 
$H_\ast(M;\mathbb{Z})$.
\end{theorem}

\noindent
Note that the Morse Homology Theorem implies that the homology of
$(C_\ast(f),\partial_\ast)$ is independent of the Morse-Smale
function $f:M \rightarrow \mathbb{R}$, the Riemannian metric,
and the chosen orientations.


\section{Morse-Bott functions and cascades}\label{cascadesection}

Let $f:M \rightarrow \mathbb{R}$ be a smooth function whose critical set
$\text{Cr}(f)$ contains a submanifold $C$ of positive dimension. 
Pick a Riemannian metric on $M$ and use it to split $T_\ast M|_C$
as
$$
T_\ast M|_C = T_\ast C \oplus \nu_\ast C
$$
where $T_\ast C$ is the tangent space of $C$ and $\nu_\ast C$ is the normal bundle of $C$.
Let $p \in C$, $V \in T_p C$, $W \in T_pM$, and let \index{Hessian} $H_p(f)$ 
be the Hessian of $f$ at $p$.  We have
$$
H_p(f)(V,W) = V_p \cdot (\tilde{W} \cdot f) = 0
$$
since $V_p \in T_pC$ and any extension of $W$ to a vector field $\tilde{W}$ satisfies 
$df(\tilde{W})|_C$ $= 0$.  Therefore, the Hessian $H_p(f)$ induces a
symmetric bilinear form $H_p^\nu(f)$ \index{$H_p^\nu(f)$ Hessian normal} on
$\nu_p C$.

\begin{definition}
A smooth function $f:M \rightarrow \mathbb{R}$ on a smooth manifold $M$ is 
called a \textbf{Morse-Bott function} if and only if the set of critical points
$\text{\rm Cr}(f)$ is a disjoint union of connected submanifolds and for each connected
submanifold $C \subseteq \text{\rm Cr}(f)$ the bilinear form $H_p^\nu(f)$ is
non-degenerate for all $p \in C$.  
\end{definition}

\noindent
Often one says that the Hessian of a Morse-Bott function $f$ is 
non-degenerate in the direction normal to the critical submanifolds.

\medskip

For a proof of the following lemma see Section 3.5 of \cite{BanLec} or
\cite{BanApr}.

\begin{lemma}[Morse-Bott Lemma] \index{Morse-Bott Lemma}\label{MorseBottLemma}
Let $f:M \rightarrow \mathbb{R}$ be a Morse-Bott function and $C \subseteq
\text{\rm Cr}(f)$ a connected component.  For any $p \in C$ there is a local chart 
of $M$ around $p$ and a local splitting $\nu_\ast C = \nu_\ast^-C \oplus 
\nu_\ast^+C$, identifying a point $x\in M$ in its domain to $(u,v,w)$ where 
$u \in C$, $v \in \nu_\ast^-C$, $w \in \nu_\ast^+C$, such that within this chart
$f$ assumes the form
$$
f(x) = f(u,v,w) = f(C)- |v|^2 + |w|^2.
$$
\end{lemma}

\begin{definition}
Let $f:M \rightarrow \mathbb{R}$ be a Morse-Bott function on a finite 
dimensional smooth manifold $M$, and let $C$ be a critical
submanifold of $f$.  For any $p \in C$ let $\lambda_p$ denote the index of 
$H_p^\nu(f)$.  This integer is the dimension of $\nu_p^-C$ 
and is locally constant by the preceding lemma.  If $C$ is connected, then 
$\lambda_p$ is constant throughout $C$ and we call $\lambda_p = \lambda_C$ 
\index{$\lambda_C$ index of $C$} the \textbf{Morse-Bott index} \index{Morse-Bott index} 
of $C$.
\end{definition}


\subsection*{Cascades}

Let $f:M \rightarrow \mathbb{R}$ be a Morse-Bott function on a finite dimensional
compact smooth manifold, and let
$$
\text{Cr}(f) = \coprod_{j=1}^l C_j,
$$
where $C_1,\ldots ,C_l$ are disjoint connected critical submanifolds of Morse-Bott index
$\lambda_1,\ldots ,\lambda_l$ respectively. Let $f_j:C_j \rightarrow
\mathbb{R}$ be a Morse function on the critical submanifold $C_j$ for all $j=1,\ldots ,l$.
If $q \in C_j$ is a critical point of $f_j:C_j\rightarrow \mathbb{R}$, then
we will denote the Morse index of $q$ relative to $f_j$ by $\lambda_q^j$, 
the stable manifold of $q$ relative to $f_j$ by $W^s_{f_j}(q) \subseteq C_j$, 
and the unstable manifold of $q$ relative to $f_j$ by $W^u_{f_j}(q) \subseteq C_j$.

\begin{definition}\label{index}
If $q\in C_j$ is a critical point of the Morse function $f_j:C_j \rightarrow \mathbb{R}$
for some $j=1,\ldots ,l$, then the \textbf{total index} of $q$, denoted $\lambda_q$, is defined 
to be the sum of the Morse-Bott index of $C_j$ and the Morse index of $q$ relative
to $f_j$, i.e.
$$
\lambda_q = \lambda_j + \lambda_q^j.
$$
\end{definition}

\noindent
The following is a restatement of Definition A.5 of \cite{FraTheA}.

\begin{definition}\label{flowlinecascade}
For $q\in \text{Cr}(f_j)$, $p \in \text{Cr}(f_i)$, and $n \in \mathbb{N}$, a 
\textbf{flow line with $n$ cascades from $q$ to $p$} is a $2n-1$-tuple:
$$
\left((x_k)_{1\leq k\leq n},(t_k)_{1 \leq k \leq n-1}  \right)
$$
where $x_k \in C^\infty(\mathbb{R},M)$ and $t_k \in \mathbb{R}_+ =
\{t \in \mathbb{R} |\ t \geq 0\}$ satisfy the following for all $k$.
\begin{enumerate}
\item Each $x_k$ is a non-constant gradient flow line of $f$, i.e.
$$
\frac{d}{dt}x_k(t) = - (\nabla f)(x_k(t)).
$$
\item For the first cascade $x_1(t)$ we have 
$$
\lim_{t \rightarrow -\infty} x_1(t) \in W^u_{f_j} (q) \subseteq C_j,
$$
and for the last cascade $x_n(t)$ we have 
$$
\lim_{t \rightarrow \infty} x_n(t) \in W^s_{f_i}(p) \subseteq C_i.
$$
\item For $1 \leq k \leq n-1$ there are critical submanifolds $C_{j_k}$ and
gradient flow lines $y_k \in C^\infty(\mathbb{R},C_{j_k})$ of $f_{j_k}$, i.e. 
$$
\frac{d}{dt}y_k(t) = - (\nabla f_{j_k})(y_k(t)),
$$
such that $\lim_{t \rightarrow \infty} x_k(t) = y_k(0)$ and 
$\lim_{t \rightarrow -\infty} x_{k+1}(t) = y_{k}(t_k)$.
\end{enumerate}
When $j = i$ a \textbf{flow line with zero cascades from $q$ to $p$} is a gradient flow
line of $f_j$ from $q$ to $p$.

\end{definition}

\begin{figure}
\includegraphics{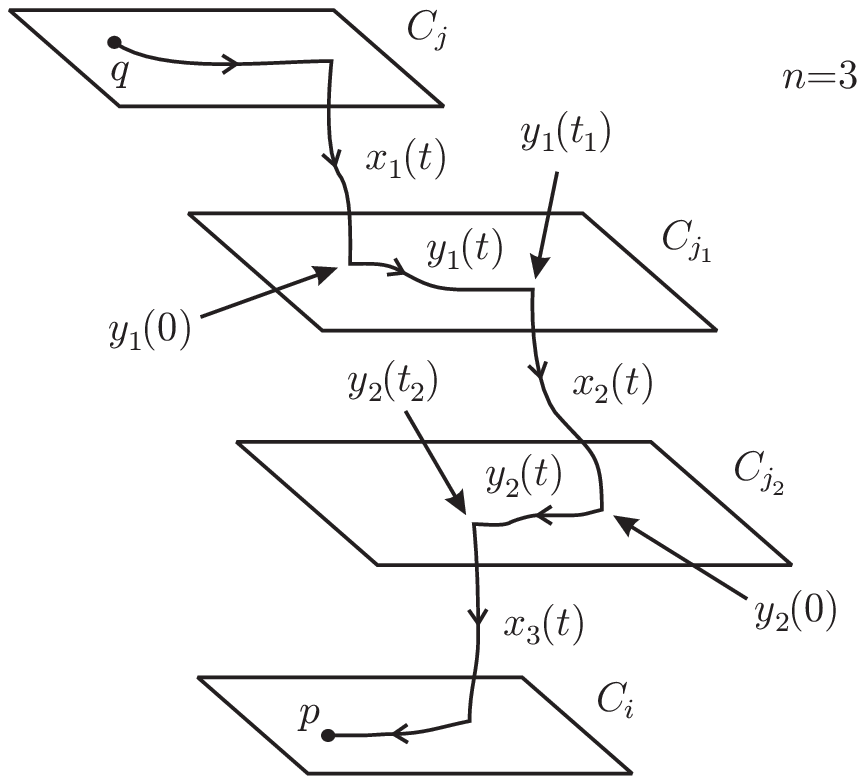}
\end{figure}

\smallskip\noindent
Note:  When $j \neq i$ a flow line with cascades from $q$ to $p$ must have at least one
cascade.

\smallskip\noindent
Note: With respect to the notation in the preceding definition, we will say that
the flow line with $n$ cascades $\left((x_k)_{1\leq k\leq n},(t_k)_{1 \leq k \leq n-1}
\right)$ \textbf{begins} at $q$ and $\textbf{ends}$ at $p$ if the conditions listed
in $(2)$ hold, i.e.
$$
\lim_{t \rightarrow -\infty} x_1(t) \in W^u_{f_j} (q) \subseteq C_j,
$$
and
$$
\lim_{t \rightarrow \infty} x_n(t) \in W^s_{f_i}(p) \subseteq C_i.
$$

\smallskip\noindent
Note: In the preceding definition the parameterizations of the gradient
flow lines $y_k(t)$ of the Morse functions $f_{j_k}:C_{j_k}\rightarrow \mathbb{R}$ are
fixed in $(3)$ by $\lim_{t \rightarrow \infty} x_k(t) = y_k(0)$,
and the entry $t_k$ records the time spent flowing along the critical submanifold
$C_{j_k}$ (or resting at a critical point). However, the parameterizations of the 
\textbf{cascades} $x_1(t),\ldots ,x_n(t)$ are not fixed. Hence, there is an action of
$\mathbb{R}^n$ on a flow line with $n$ cascades given by
$$
\left((x_k(t))_{1\leq k\leq n},(t_k)_{1 \leq k \leq n-1}  \right) \mapsto
\left((x_k(t+s_k))_{1\leq k\leq n},(t_k)_{1 \leq k \leq n-1}  \right)
$$
for $(s_1,\ldots ,s_n) \in \mathbb{R}^n$.

\begin{definition}
For $q\in \text{Cr}(f_j)$, $p \in \text{Cr}(f_i)$, and $n \in \mathbb{N}$ we denote
the space of flow lines from $q$ to $p$ with $n$ cascades by $W^c_n(q,p)$, and we 
denote the quotient of $W^c_n(q,p)$ by the action of $\mathbb{R}^n$ by
$$
\mathcal{M}^c_n(q,p) = W^c_n(q,p)/\mathbb{R}^n.
$$
The \textbf{set of unparameterized flow lines with cascades from $q$ to $p$} is defined
to be
$$
\mathcal{M}^c(q,p) = \bigcup_{n \in \mathbb{Z}_+} \mathcal{M}^c_n(q,p)
$$
where $\mathcal{M}^c_0(q,p) = W^c_0(q,p)/\mathbb{R}$. We will say 
that an element of $\mathcal{M}^c(q,p)$ \textbf{begins} at $q$ and \textbf{ends} at $p$.
\end{definition}

\smallskip
We now prove that $\mathcal{M}^c(q,p)$ is a smooth manifold of dimension 
$\lambda_q - \lambda_p - 1$ when $f:M \rightarrow \mathbb{R}$ satisfies the Morse-Bott-Smale
transversality condition with respect to the metric $g$, the Morse functions $f_k:C_k
\rightarrow \mathbb{R}$ satisfy the Morse-Smale transversality condition with respect to
the restriction of $g$ to $C_k$ for all $k=1,\ldots ,l$, and the stable and unstable
manifolds of the Morse-Smale functions $f_i:C_i \rightarrow \mathbb{R}$ and 
$f_j:C_j \rightarrow \mathbb{R}$ are transverse to certain beginning and endpoint maps.
Our proof uses fibered product constructions on smooth manifolds with corners similar to those
found in \cite{BanMor}.

\begin{definition}[Morse-Bott-Smale Transversality]\label{MBStransversality}
A Morse-Bott function $f:M \rightarrow \mathbb{R}$ is said to satisfy the 
\textbf{Morse-Bott-Smale transversality} condition with respect to a given 
Riemannian metric $g$ on $M$ if and only if for any two connected critical submanifolds 
$C$ and $C'$, $W^u_f(q)$ intersects $W^s_f(C')$ transversely in $M$, i.e. $W^u_f(q)
\pitchfork W^s_f(C') \subseteq M$, for all $q\in C$.
\end{definition}

Let $C_k$ and $C_{k'}$ be two connected critical submanifolds of $f$, and let $W^u_f(C_k)$
and $W^s_f(C_{k'})$ denote the unstable and stable manifolds of $C_k$ and $C_{k'}$
with respect to the flow of $-\nabla f$. The Morse-Bott-Smale transversality assumption 
implies that the moduli space of gradient flow lines of $f$:
$$
\mathcal{M}_f(C_k,C_{k'}) = \left(W^u_f(C_k) \cap W^s_f(C_{k'}) \right) /\mathbb{R}
$$
is either empty or a smooth manifold of dimension $\lambda_k - \lambda_{k'} + 
\text{dim }C_k - 1$. Moreover, the beginning and endpoint maps 
$\partial_-:\mathcal{M}_f(C_k,C_{k'}) \rightarrow C_k$ and 
$\partial_+:\mathcal{M}_f(C_k,C_{k'}) \rightarrow C_{k'}$ are smooth, and the beginning
point map $\partial_-$ is a submersion (see Lemma 5.19 of \cite{BanMor}). 

Now assume that the following moduli spaces and fibered products are nonempty.
Then for distinct $k,k',k'' \in \{1,2,\ldots ,l\}$ and $t \in \mathbb{R}_+ = 
\{t \in \mathbb{R} |\ t \geq 0\}$ we can consider the fibered product
$$
\xymatrix{
(\mathbb{R}_+ \times \mathcal{M}_f(C_{k},C_{k'}) ) \times_{C_{k'}} 
\mathcal{M}_f(C_{k'},C_{k''}) \ar@{-->}[r] \ar@{-->}[d] & 
\mathcal{M}_f(C_{k'},C_{k''}) \ar[d]^{\partial_-}\\
\mathbb{R}_+ \times \mathcal{M}_f(C_{k},C_{k'}) 
\ar[r]^(.68){\varphi_t \circ \partial_+ \circ \pi_2} & C_{k'}
}
$$
where $\pi_2$ denotes projection onto the second component and $\varphi_t$ denotes the gradient 
flow of $f_{k'}$ along the critical submanifold $C_{k'}$ for time $t \in \mathbb{R}_+$.
This fibered product is a smooth manifold with boundary because $\partial_-:\mathcal{M}_f
(C_{k'},C_{k''}) \rightarrow C_{k'}$ is a submersion, and its dimension is
\begin{eqnarray*}
& & (\lambda_k - \lambda_{k'} + \text{dim }C_k) + (\lambda_{k'} - \lambda_{k''} + 
\text{dim }C_{k'} - 1) - \text{dim }C_{k'}\\
& = & \lambda_k - \lambda_{k''} + \text{dim }C_{k} - 1
\end{eqnarray*}
(see Lemma 4.5 and Lemma 5.21 of \cite{BanMor}). Similarly, for any set of distinct
integers $\{j_1, j_2,\ldots ,j_{n-1}\} \subseteq \{1,2,\ldots ,l\}$ such that the
following moduli spaces are nonempty, the iterated fibered product
\begin{eqnarray*}
(\mathbb{R}_+ \times \mathcal{M}_f(C_j,C_{j_1})) \times_{C_{j_1}} 
(\mathbb{R}_+ \times \mathcal{M}_f(C_{j_1},C_{j_2})) \times_{C_{j_2}} \cdots \\
\times_{C_{j_{n-2}}} (\mathbb{R}_+ \times \mathcal{M}_f(C_{j_{n-2}},C_{j_{n-1}})) 
\times_{C_{j_{n-1}}} \mathcal{M}_f(C_{j_{n-1}},C_i)
\end{eqnarray*}
is a smooth manifold with corners because $\partial_-\circ \pi_2: \mathbb{R}_+ \times
\mathcal{M}_f (C_{k},C_{k'}) \rightarrow C_k$ is a submersion and a stratum submersion
for all $k,k'=1,\ldots ,l$. We will denote this smooth manifold with corners by 
$\mathcal{M}^c_n(C_j,C_{j_1},\ldots ,C_{j_{n-1}},C_i)$. Its dimension is
\begin{eqnarray*}
&   & (\lambda_j - \lambda_{j_1} + \text{dim }C_j) + (\lambda_{j_1} - \lambda_{j_2} 
      + \text{dim }C_{j_1}) - \text{dim }C_{j_1} + \cdots\\
& + & (\lambda_{j_{n-2}} - \lambda_{j_{n-1}} + \text{dim }C_{j_{n-2}})
       - \text{dim }C_{j_{n-2}}+ (\lambda_{j_{n-1}} - \lambda_i + \text{dim }C_{j_{n-1}} - 1)\\ 
& - & \text{dim }C_{{j_{n-1}}}\\
= & & \lambda_j - \lambda_i + \text{dim }C_j - 1,
\end{eqnarray*}
which is independent of $j_1,j_2,\ldots ,j_{n-1}$.  Note that we have smooth
beginning and endpoint maps
\begin{eqnarray*}
\partial_-:\mathcal{M}^c_n(C_j,C_{j_1},\ldots ,C_{j_{n-1}},C_i) \rightarrow C_j\\
\partial_+:\mathcal{M}^c_n(C_j,C_{j_1},\ldots ,C_{j_{n-1}},C_i) \rightarrow C_i.
\end{eqnarray*}

\smallskip
We can now state our transversality assumptions for the stable and unstable manifolds
$W^s_{f_{i}}(p)$ and $W^u_{f_j}(q)$ of the Morse-Smale functions 
$f_i:C_i \rightarrow \mathbb{R}$ and $f_j:C_j \rightarrow \mathbb{R}$ with respect 
to these beginning and endpoint maps.

\begin{definition}\label{beginningendtransverse}
The stable and unstable manifolds $W^s_{f_{i}}(p)$ and $W^u_{f_j}(q)$ are
\textbf{transverse to the beginning and endpoint maps} if and only if for any set
(possibly empty) of distinct integers $\{j_1, j_2,\ldots ,j_{n-1}\} \subseteq
\{1,2,\ldots ,l\}$ such that the moduli space $\mathcal{M}^c_n(C_j,C_{j_1},\ldots ,C_{j_{n-1}},
C_i)$ is not empty the map
$$
\mathcal{M}^c_n(C_j,C_{j_1},\ldots ,C_{j_{n-1}},C_i) \stackrel{(\partial_-,\partial_+)}
{\longrightarrow} C_j \times C_{i}
$$
is transverse and stratum transverse to $W^u_{f_j} (q) \times W^s_{f_{i}}(p)$.
\end{definition}

\smallskip\noindent
Note:  When $\{j_1, j_2,\ldots ,j_{n-1}\} = \emptyset$ we have $\mathcal{M}^c_1(C_j,C_i)
 = \mathcal{M}_f(C_j,C_i)$.

\begin{lemma}\label{perturb}
There exist arbitrarily small perturbations of $f_i:C_i \rightarrow \mathbb{R}$ and
$f_j:C_j \rightarrow \mathbb{R}$ to smooth Morse-Smale functions $\tilde{f}_i$ and 
$\tilde{f}_j$ such that all the stable and unstable manifolds of $\tilde{f}_i$ and 
$\tilde{f}_j$ are transverse to the beginning and endpoint maps. Moreover, there exist
open neighborhoods of $\tilde{f}_i$ and $\tilde{f}_j$ consisting of smooth Morse-Smale 
functions whose stable and unstable manifolds are all transverse to the beginning and 
endpoint maps.
\end{lemma}

\proofstart
Let $\{j_1, j_2,\ldots ,j_{n-1}\} \subseteq \{1,2,\ldots ,l\}$ be a (possibly empty) set of
distinct integers such that the moduli space $\mathcal{M}^c_n(C_j,C_{j_1},\ldots ,C_{j_{n-1}},
C_i)$ is not empty, and let $X$ be a stratum of $\mathcal{M}^c_n(C_j,C_{j_1},\ldots ,
C_{j_{n-1}}, C_i)$. Let
$$
E^s_{f_i}:\mathbb{R}^{\text{dim }C_i - \lambda_p^i} \rightarrow W^s_{f_i}(p) \subseteq C_i
$$
and
$$
E^u_{f_j}:\mathbb{R}^{\lambda_q^j} \rightarrow W^u_{f_j}(q)\subseteq C_j
$$
be the surjective smooth embeddings from Section \ref{MorseChain}, where $p \in \text{Cr}(f_i)$,
$q \in \text{Cr}(f_j)$, and we have identified $T^s_pC_i = \mathbb{R}^{\text{dim }C_i - 
\lambda_p^i}$ and $T^u_qC_j = \mathbb{R}^{\lambda_q^j}$. The stable and unstable manifolds
$W^s_{f_j}(p)$ and $W^u_{f_i}(q)$ are transverse to $(\partial_-,\partial_+):X 
\rightarrow C_j \times C_i$ if and only if the map
$$
(E^u_{f_j}, E^s_{f_i}) \times (\partial_-,\partial_+): (\mathbb{R}^{\lambda_q^j}
\times \mathbb{R}^{\text{dim }C_i - \lambda_p^i}) \times X \rightarrow (C_j \times C_i) \times 
(C_j \times C_i)
$$
is transverse to the diagonal $\Delta \subset (C_j \times C_i) \times (C_j \times C_i)$.

For any $r \geq 2$ the set of $C^r$ Morse-Smale functions on a smooth Riemannian manifold 
$(M,g)$ is an open and dense subset of the set of all $C^r$ functions on $M$, and the phase
diagram of a Morse-Smale function is stable under small $C^r$ perturbations \cite{PalOnM}. Thus,
there exist a neighborhood $\mathcal{N}_{f_i}\subset C^r(M,\mathbb{R})$ of $f_i$ 
such that $\tilde{f}_i \in \mathcal{N}_{f_i}$ implies that $\tilde{f}_i$ is a Morse-Smale
function with critical points of the same index and near the critical points of $f_i$.
Similarly, there exists a neighborhood $\mathcal{N}_{f_j} \subset C^r(M,\mathbb{R})$ of $f_j$
such that $\tilde{f}_j \in \mathcal{N}_{f_j}$ implies that $\tilde{f}_j$ is a Morse-Smale
function with critical points of the same index and near the critical points of $f_i$.
Moreover, we can choose these neighborhoods small enough so that the maps
$$
E^s:\mathcal{N}_{f_i} \rightarrow C^r(\mathbb{R}^{\text{dim }C_i - \lambda_p^i},C_i)
\quad \text{ and } \quad E^u:\mathcal{N}_{f_j} \rightarrow C^r(\mathbb{R}^{\lambda_q^j}, C_j)
$$
defined by sending $\tilde{f}_i \in \mathcal{N}_{f_i}$ to the embedding $E^s_{\tilde{f}_i}$
(with respect to the critical point $\tilde{p}$ near $p$) and $\tilde{f}_j \in \mathcal{N}_{f_j}$
to the embedding $E^u_{\tilde{f}_j}$ (with respect to the critical point $\tilde{q}$ near $q$)
are well defined and of class $C^r$. In particular, we can choose the neighborhoods small enough 
so that we can identify $T^s_{\tilde{p}}C_i = T^s_pC_i = \mathbb{R}^{\text{dim }C_i - \lambda_p^i}$
and $T^u_{\tilde{q}}C_j = T^u_qC_j = \mathbb{R}^{\lambda_q^j}$.

\smallskip
The map
$$
(E^u \times E^s) \times (\partial_-,\partial_+): (\mathcal{N}_{f_j} \times \mathcal{N}_{f_i})
\times (\mathbb{R}^{\lambda_q^j} \times \mathbb{R}^{\text{dim }C_i - \lambda_p^i} \times X) 
\rightarrow (C_j \times C_i) \times (C_j \times C_i)
$$
defined by
$$
\left((E^u \times E^s) \times (\partial_-,\partial_+)\right) \left((\tilde{f}_j,\tilde{f}_i)
\times (x,y,\gamma) \right) = (E^u_{\tilde{f}_j}(x), E^s_{\tilde{f}_i}(y)) \times
(\partial_-(\gamma), \partial_+(\gamma))
$$
is of class $C^r$ (see Theorem 12.3 of \cite{AbrTra}) and transverse to $\Delta \subset 
(C_j \times C_i) \times (C_j \times C_i)$.  Hence, by the Transversality Density Theorem
(Theorem 19.1 of \cite{AbrTra}) the set of Morse-Smale functions $(\tilde{f}_j,\tilde{f_i}) 
\in \mathcal{N}_{f_j} \times \mathcal{N}_{f_i}$ such that  
$$
(E^u_{\tilde{f}_j}, E^s_{\tilde{f}_i}) \times (\partial_-, \partial_+): (\mathbb{R}^{\lambda_q^j}
\times \mathbb{R}^{\text{dim }C_i - \lambda_p^i}) \times X \rightarrow
(C_j \times C_i) \times (C_j \times C_i)
$$
is transverse to $\Delta$ is residual (and hence dense) in $\mathcal{N}_{f_j} \times 
\mathcal{N}_{f_i}$ for $r \geq 2$ large enough, e.g. $r > 3\ \text{dim }M$.  

Since there are only finitely many subsets $\{j_1, j_2,\ldots ,j_{n-1}\} \subseteq 
\{1,2,\ldots ,l\}$, finitely many critical points of $f_i$ and $f_j$, and finitely
many strata $X$, we can intersect finitely many such residual sets to obtain a residual
(and hence dense) subset $\mathcal{R} \subseteq\mathcal{N}_{f_j} \times
\mathcal{N}_{f_i}$ such that 
$$
\mathcal{M}^c_n(C_j,C_{j_1},\ldots ,C_{j_{n-1}},C_i) \stackrel{(\partial_-,\partial_+)}
{\longrightarrow} C_j \times C_{i}
$$
is transverse and stratum transverse to $W^u_{\tilde{f}_j}(\tilde{q}) \times 
W^s_{\tilde{f}_i}(\tilde{p})$ for all $\tilde{q} \in \text{Cr}(\tilde{f}_j)$ and
$\tilde{p} \in \text{Cr}(\tilde{f}_i)$ whenever $(\tilde{f}_j,\tilde{f}_i) \in \mathcal{R}$.
Moreover, since the space of smooth Morse-Smale functions on $M$ is dense in the space of
$C^r$ Morse-Smale functions on $M$, the Openness of Transversal Intersection Theorem
(Theorem 18.2 of \cite{AbrTra}) implies that we can find open neighborhoods of smooth
functions arbitrarily close to $f_j$ and $f_i$ consisting of Morse-Smale functions 
$\tilde{f}_j$ and $\tilde{f}_i$ with $(\tilde{f}_j,\tilde{f}_i) \in \mathcal{R}$.

\proofend

\smallskip\noindent
Note: The critical points of $f_i$ and $f_j$ may not be preserved by the perturbations
in the preceding lemma. However, it is possible to choose the perturbations so that
the phase diagrams of $f_i$ and $f_j$ do not change \cite{PalOnM}. In particular, the number
of critical points of index $k$ remains the same for all $k=1, \ldots, m$, which also
follows from the Rigidity Theorem 1.19 of \cite{CorRig}. 

\smallskip
The next theorem should be compared with Theorem A.12 from \cite{FraTheA}, whose proof
uses the modern infinite dimensional techniques of Floer homology. Theorem A.12 in
\cite{FraTheA} is proved under the assumption that the Riemannian metric $g$ on $M$ is
generic, which is necessary to ensure that a certain Fredholm operator used in the proof
of the theorem is surjective. 


\begin{theorem}\label{manifold}
Assume that $f$ satisfies the Morse-Bott-Smale transversality condition with respect to
the Riemannian metric $g$ on $M$, $f_k:C_k \rightarrow \mathbb{R}$ satisfies
the Morse-Smale transversality condition with respect to the restriction of $g$
to $C_k$ for all $k=1,\ldots ,l$, and the unstable and stable manifolds $W^u_{f_j}(q)$ and
$W^s_{f_{i}}(p)$ are transverse to the beginning and endpoint maps. 
\begin{enumerate}
\item When $n=0,1$ the set $\mathcal{M}^c_n(q,p)$ is either empty or a smooth manifold
      without boundary.
\item For $n > 1$ the set $\mathcal{M}^c_n(q,p)$ is either empty or a smooth manifold with
      corners.
\item The set $\mathcal{M}^c(q,p)$ is either empty or a smooth manifold without boundary.
\end{enumerate}
In each case the dimension of the manifold is $\lambda_q - \lambda_p - 1$.
When $M$ is orientable and $C_k$ is orientable for all $k=1,\ldots, l$, the above manifolds
are orientable.
\end{theorem}

\proofstart
For more details concerning the notation and dimension formulas used in the following
we refer the reader to Sections 3 and 4 of \cite{BanMor}. We first prove statements (1) and
(2) using pullback constructions. A gluing theorem is then used to show that the space
$\mathcal{M}^c_{\leq n}(C_j,C_i)$ consisting of flow lines with at most $n$ cascades beginning
at any point in $C_j$ and ending at any point in $C_i$ is a manifold without boundary.
Pulling back $W^u_{f_j}(q) \times W^s_{f_{i}}(p)$ via the beginning and endpoint maps on
$\mathcal{M}^c_{\leq l}(C_j,C_i)$ then shows that $\mathcal{M}^c(q,p)$ is a smooth manifold
without boundary of dimension $\lambda_q - \lambda_p - 1$.

\smallskip
The space $\mathcal{M}_0^c(q,p)$ is empty unless $i=j$, and when $i=j$ the theorem 
follows from the fact that $f_j$ satisfies the Morse-Smale transversality condition. 
For the case $n=1$ note that the assumption that
$$
\mathcal{M}^c_1(C_j,C_i) \stackrel{(\partial_-,\partial_+)}{\longrightarrow} C_j \times C_{i}
$$
is transverse to $W^u_{f_j}(q) \times W^s_{f_{i}}(p)$ implies that
$$
\mathcal{M}^c_1(W^u_{f_j}(q),W^s_{f_i}(p)) \stackrel{\text{def}}{=} 
(\partial_-,\partial_+)^{-1} (W^u_{f_j}(q) \times W^s_{f_i}(p))
$$
is either empty or a smooth manifold. In the second case, the codimension of the
manifold $\mathcal{M}^c_1(W^u_{f_j}(q),W^s_{f_i}(p))$ is $\text{dim } C_j - \lambda^j_q + 
\lambda_p^i$, and hence the dimension of $\mathcal{M}^c_1(W^u_{f_j}(q),W^s_{f_i}(p))$ is 
$\lambda_j + \lambda^j_q - (\lambda_i + \lambda^i_p) - 1$ since the dimension of 
$\mathcal{M}^c_1(C_j,C_i)$ is $\lambda_j - \lambda_i + \text{dim } C_j - 1$. 
(See for instance Theorem 5.11 of \cite{BanLec}.) This shows that 
$\mathcal{M}^c_1(q,p) = \mathcal{M}^c_1(W^u_{f_j}(q),W^s_{f_i}(p))$ is a smooth manifold
without boundary of dimension $\lambda_q - \lambda_p - 1$.

\smallskip
Now assume that $n > 1$ and the following moduli spaces and fibered products
are nonempty. Then for distinct $j_1, j_2,\ldots ,j_{n-1} \in \{1,2,\ldots ,l\}$
the assumption that
$$
\mathcal{M}^c_n(C_j,C_{j_1},\ldots ,C_{j_{n-1}}, C_i) \stackrel{(\partial_-,\partial_+)}
{\longrightarrow} C_j \times C_{i}
$$
is transverse and stratum transverse to $W^u_{f_j} (q) \times W^s_{f_{i}}(p)$ implies that
$$
\mathcal{M}^c_n(W^u_{f_j}(q), C_{j_1},\ldots , C_{j_{n-1}},W^s_{f_i}(p)) 
\stackrel{\text{def}}{=} (\partial_-,\partial_+)^{-1} (W^u_{f_j}(q) \times W^s_{f_i}(p))
$$
is a smooth manifold with corners of dimension $\lambda_q - \lambda_p - 1$. This shows that
$$
\mathcal{M}^c_n(q,p)\ = \bigcup_{\{j_1,\ldots ,j_{n-1}\}} \mathcal{M}^c_n(W^u_{f_j}(q),
C_{j_1},\ldots , C_{j_{n-1}},W^s_{f_i}(p))
$$
is a smooth manifold with corners of dimension $\lambda_q - \lambda_p - 1$, where the union 
is taken over all sets of distinct integers $\{j_1,\ldots ,j_{n-1}\} \subseteq 
\{1,2,\ldots ,l\}$. This completes the proof of statements (1) and (2).

\smallskip
We now use a gluing theorem to define smooth charts on 
$$
\mathcal{M}_{\leq n}^c(C_j,C_i) \stackrel{\text{def}}{=} \bigcup_{k=0}^n 
\mathcal{M}_{k}^c(C_j,C_i)
$$
where $\mathcal{M}_{k}^c(C_j,C_i)$ denotes the union of $\mathcal{M}^c_k(C_j, C_{j_1},
\ldots, C_{j_{k-1}},C_i)$ over all sets of distinct integers $\{j_1,\ldots ,j_{k-1}\} 
\subseteq \{1,2,\ldots ,l\}$ when $k > 1$.  For distinct $k,k',k'' \in \{1,2,\ldots , l\}$
there exists an $\varepsilon >0$ and a smooth injective local diffeomorphism
$$
G:\mathcal{M}_f(C_k,C_{k'}) \times_{C_{k'}} \mathcal{M}_f(C_{k'},C_{k''}) \times 
(-\varepsilon,0) \rightarrow \mathcal{M}_f(C_k,C_{k''})
$$
onto an end of $\mathcal{M}_f(C_k,C_{k''})$, where the fibered product is taken with respect
to the beginning and endpoint maps $\partial_-$ and $\partial_+$. (See for instance 
Appendix A.3 of \cite{AusMor} or Theorem 4.8 of \cite{BanMor}.) 
Let $\rho:(-\varepsilon,\infty) \rightarrow (-\varepsilon,\infty)$ be a smooth
map that is smoothly homotopic to
$$
\chi(t) = \left\{
\begin{array}{ll}
t & t \geq 0 \\
0 & t \leq 0
\end{array}
\right.
$$
and satisfies 
$$
\rho(t) = \left\{
\begin{array}{ll}
t & t \geq \varepsilon/2 \\
0 & t \leq 0.
\end{array}
\right.
$$
For $\varepsilon > 0$ sufficiently small we can replace the maps 
$\varphi_t \circ \partial_+ \circ \pi_2$ in the iterated fibered product that defines
$\mathcal{M}^c_n(C_j,C_{j_1},\ldots ,C_{j_{n-1}}, C_i)$ with the maps $\varphi_{\rho(t)}
\circ \partial_+ \circ \pi_2$ and obtain a smooth manifold with corners that is smoothly
diffeomorphic to the original manifold. Moreover, if we choose $\varepsilon > 0$ small enough,
then $W^u_{f_j}(q)$ and $W^s_{f_{i}}(p)$ will still be transverse to the beginning and
endpoint maps from the modified fibered product space.

\smallskip
Using the maps $\varphi_{\rho(t)} \circ \partial_+ \circ \pi_2$ and $\partial_-$ we
consider the fibered product 
$$
((-\varepsilon,\infty) \times \mathcal{M}_f(C_{j},C_{k}) ) \times_{C_{k}} 
\mathcal{M}_f(C_{k},C_{i})
$$
where $k \in \{1,2,\ldots ,l\}$. The part of this smooth manifold where
$-\varepsilon < t < 0$ is diffeomorphic to an end of $\mathcal{M}^c_1(C_j,C_i)$ by the
above gluing theorem, and the part of the space where $t \geq 0$ is diffeomorphic to
$\mathcal{M}^c_2(C_j,C_k,C_i)$. Therefore, there are smooth charts on the above manifold
around the points where $t=0$ which are compatible with the smooth charts
on $\mathcal{M}^c_1(C_j,C_i)$ and the smooth charts on $\mathcal{M}^c_2(C_j,C_k,C_i)$.
This shows that the space $\mathcal{M}^c_{\leq 2}(C_j,C_i)$ of unparameterized flow lines
with at most $2$ cascades from $C_j$ to $C_i$ is a smooth manifold without boundary of
dimension $\lambda_j - \lambda_i + \text{dim } C_j - 1$.

Continuing by induction, for distinct $j_1,j_2,\ldots ,j_{n-1} \in \{1,2,\ldots ,l\}$
the fibered product
\begin{eqnarray*}
((-\varepsilon,\infty) \times \mathcal{M}_f(C_j,C_{j_1})) \times_{C_{j_1}} 
((-\varepsilon,\infty) \times \mathcal{M}_f(C_{j_1},C_{j_2})) \times_{C_{j_2}} \cdots \\
\times_{C_{j_{n-2}}} ((-\varepsilon,\infty) \times \mathcal{M}_f(C_{j_{n-2}},C_{j_{n-1}})) 
\times_{C_{j_{n-1}}} \mathcal{M}_f(C_{j_{n-1}},C_i)
\end{eqnarray*}
with respect to the maps $\varphi_{\rho(t)} \circ \partial_+ \circ \pi_2$ and 
$\partial_- \circ \pi_2$ is a smooth manifold.  The part of the space where 
$-\varepsilon < t_k < 0$ for some $k$ is diffeomorphic to an end of $\mathcal{M}^c_{\leq n-1}
(C_j,C_i)$, and the part of the space where $t_k \geq 0$ for all $k$ is diffeomorphic to
$\mathcal{M}^c_n(C_j, C_{j_1},\ldots , C_{j_{n-1}},C_i)$. Thus, the space
$\mathcal{M}^c_{\leq n}(C_j,C_i)$ of unparameterized flow lines with at most $n$ cascades from 
$C_j$ to $C_i$ is a smooth manifold without boundary of dimension $\lambda_j - \lambda_i + 
\text{dim } C_j - 1$.  Moreover, 
$$
\mathcal{M}^c_{\leq n}(C_j,C_i) \stackrel{(\partial_-,\partial_+)}{\longrightarrow}
C_j \times C_{i}
$$
is transverse to $W^u_{f_j} (q) \times W^s_{f_{i}}(p)$. The pullback of 
$W^u_{f_j} (q) \times W^s_{f_{i}}(p)$ under this map is the space of
unparameterized flow lines with at most $n$ cascades from $q$ to $p$:
$$
\mathcal{M}^c_{\leq n}(q,p) = \bigcup_{k=0}^n \mathcal{M}^c_k(q,p).
$$
Hence, for any $0 \leq n \leq l$ the space $\mathcal{M}^c_{\leq n}(q,p)$ is either
empty or a smooth manifold without boundary of dimension $\lambda_q - \lambda_p - 1$. 
Taking $n=l$ we see that $\mathcal{M}^c(q,p)$ is either empty or a smooth manifold without
boundary of dimension $\lambda_q - \lambda_p - 1$.


Now, an orientation on $M$ and orientations on $C_j$ for all $j=1,\ldots ,l$ determine
orientations on the above fibered products by the results in Section 5.2 of \cite{BanMor}.
If we choose the gluing diffeomorphisms to be compatible with these orientations,
then we obtain an orientation on $\mathcal{M}^c(q,p)$.

\proofend


\section{Broken flow lines with cascades}\label{broken}

We will now consider the compactness properties of $\mathcal{M}^c(q,p)$. In general,
$\mathcal{M}^c(q,p)$ will be a non-compact manifold because a sequence of unparameterized
flow lines with cascades from $q$ to $p$ may converge to a broken flow line with cascades from
$q$ to $p$. Throughout this section we will assume that $f$ satisfies the Morse-Bott-Smale
transversality condition with respect to the Riemannian metric $g$ on a compact smooth
manifold $M$, $f_k:C_k \rightarrow \mathbb{R}$ satisfies the Morse-Smale transversality
condition with respect to the restriction of $g$ to $C_k$ for all $k=1,\ldots ,l$, and the
unstable and stable manifolds $W^u_{f_j}(q)$ and $W^s_{f_{i}}(p)$ are transverse to the
beginning and endpoint maps. 

\medskip
It is well known that any sequence of unparameterized gradient flow lines between two critical
points of a Morse-Smale function must have a subsequence that converges to a broken flow line.
However, making this statement precise requires a discussion of the topology on the space of
broken flow lines. The topology on the space of broken flow lines can be defined in several
ways, including the compact open topology (after picking specific parameterizations for the
flow lines), in terms of Floer-Gromov convergence, and using the Hausdorff metric (after
identifying a broken flow line with its image). For a detailed discussion concerning different
ways to define the topology on the space of broken flow lines of a Morse-Smale function
and proofs that the resulting spaces are homeomorphic see \cite{LetMor}.

\smallskip
To prove a similar result for cascades we first need to explain what we mean by a
broken flow line with cascades. Roughly speaking, a broken flow line with cascades 
is an unparameterized flow line with cascades that either flows along an intermediate
critical submanifold for infinite time or rests at an intermediate critical point
of one of the Morse functions $f_k:C_k \rightarrow \mathbb{R}$ for some $k=1,\ldots, l$
for infinite time. To make this more precise, recall that a flow line with
cascades is of the form $((x_k)_{1 \leq k \leq n}, (t_k)_{1 \leq k \leq n-1})$ where
$t_k \in \mathbb{R}_+ = \{t \in \mathbb{R}|\ t \geq 0\}$. In particular, $t_k < \infty$, 
but we might have $t_k = 0$ for some $k$. If $t_k = 0$ for some $k$, then the
flow line with cascades ``looks like'' it contains a broken flow line. That is, if 
$t_k = 0$, then $\lim_{t \rightarrow \infty} x_k(t)= \lim_{t \rightarrow -\infty} 
x_{k+1}(t)$ and $(x_k,x_{k+1})$ is a broken flow line of the Morse-Bott function 
$f:M \rightarrow \mathbb{R}$.  However, $(x_k,x_{k+1},0)$ is an unbroken flow line with 
$2$ cascades. 

Since a flow line with cascades must begin and end at critical points of the Morse
functions chosen on the critical submanifolds, it's clear that $(x_k,x_{k+1})$ should not
be called a broken flow line with cascades when $\lim_{t \rightarrow \infty}
x_k(t)= \lim_{t \rightarrow -\infty} x_{k+1}(t)$ is not a critical point of 
$f_{j_k}:C_{j_k} \rightarrow \mathbb{R}$. In order to be consistent, we will \textbf{not}
call $(x_k,x_{k+1})$ a broken flow line with cascades even if $\lim_{t \rightarrow \infty}
x_k(t) = \lim_{t \rightarrow -\infty} x_{k+1}(t) = r$ is a critical point of 
$f_{j_k}:C_{j_k} \rightarrow \mathbb{R}$.  Instead, we will always assume that the 
time spent resting at the intermediate critical point is zero, unless the time is 
otherwise specified. That is, we will identify $(x_k,x_{k+1})$ with the flow line
with $2$ cascades $(x_k,x_{k+1},0)$.

In general, suppose that we have  an $n$-tuple of unparameterized flow lines with cascades
$(v_1,\ldots ,v_n)$ such that $v_1$ begins at $q \in \text{Cr}(f_j)$, $v_n$ ends at
$p\in \text{Cr}(f_i)$, and $v_\nu$ begins where $v_{\nu-1}$ ends for $2\leq \nu \leq n$.
Suppose that $v_\nu$ is represented by $((x_k^\nu)_{1 \leq k \leq n_\nu}, (t_k^\nu)_{1 \leq k 
\leq n_\nu-1})$ and $v_{\nu-1}$ is represented by $((x_k^{\nu-1})_{1 \leq k \leq n_{\nu-1}},
(t_k^{\nu-1})_{1 \leq k \leq n_{\nu-1}-1})$. The statement that $v_\nu$ begins
where $v_{\nu - 1}$ ends means that there is a critical point $r$ of one of the
Morse functions $f_k:C_k \rightarrow \mathbb{R}$ for some $k = 1, \ldots, l$ such that 
$\lim_{t \rightarrow \infty} x_{n_{\nu-1}}^{\nu-1}(t) \in W^s_{f_j}(r)$ and 
$\lim_{t \rightarrow -\infty} x_{n_\nu}^\nu(t) \in W^u_{f_j}(r)$.
So, it appears that $(v_{\nu-1},v_\nu)$ differs from an unparameterized flow line with cascades
in that $(v_{\nu-1},v_\nu)$ flows along the intermediate critical submanifold $C_k$ for 
infinite time. However, if $\lim_{t \rightarrow \infty} x_{n_{\nu-1}}^{\nu-1}(t) =  
\lim_{t \rightarrow -\infty} x_{n_\nu}^\nu(t) = r$, then $(v_{\nu-1},v_\nu)$ determines an
unparameterized flow line with $n_{\nu-1} + n_\nu$ cascades where the time spent resting at
the intermediate critical point $q$ is $0$, i.e. the unparameterized flow line with cascades
represented by 
$$
((x_k^{\nu-1})_{1 \leq k \leq n_{\nu-1}},(x_k^\nu)_{1 \leq k \leq n_\nu}, 
(t_k^{\nu-1})_{1 \leq k \leq n_{\nu-1}-1},0,(t_k^\nu)_{1 \leq k \leq n_\nu-1}).
$$
In this case, we will identify $(v_{\nu-1},v_\nu)$ with the unparameterized flow line with
cascades represented by the above tuple.

\smallskip
It is interesting to consider what this convention means for a Morse-Smale function
$f:M \rightarrow \mathbb{R}$. Suppose that $p,q,r \in \text{Cr}(f)$, $\gamma_1$ is
a gradient flow line from $q$ to $r$ and $\gamma_2$ is a gradient flow line from
$r$ to $p$.  Then with this convention we are identifying the broken gradient
flow line represented by $(\gamma_1,\gamma_2)$ with the flow line with $2$ cascades
$(\gamma_1,\gamma_2,0)$. In fact, for a Morse-Smale function this convention means
that the only truly broken flow lines with cascades have representations of the form 
$((x_k)_{1 \leq k \leq n}, (t_k)_{1 \leq k \leq n-1})$, where $t_k = \infty$ for some $k$.

\begin{definition}
A \textbf{broken flow line with cascades} from $q\in Cr(f_j)$ to $p \in Cr(f_i)$ is an 
$n$-tuple of unparameterized flow lines with cascades $(v_1,\ldots ,v_n)$ such that $v_1$
begins at $q$, $v_n$ ends at $p$, and $v_\nu$ begins where $v_{\nu-1}$ ends for $2\leq \nu 
\leq n$, subject to the following restriction. If the last cascade of $v_{\nu-1}$ and the
first cascade of $v_\nu$ meet at a critical point of one of the Morse functions
$f_k:C_k \rightarrow \mathbb{R}$ for some $k=1,\ldots, l$, then the time spent resting 
at the critical point is infinity.
\end{definition}

A sequence of unparameterized flow lines with cascades from $q\in Cr(f_j)$ to $p \in Cr(f_i)$
must have a subsequence that converges to a broken flow line with cascades from $q$ to $p$.
This is proved in Theorem A.10 of \cite{FraTheA} with respect to Floer-Gromov convergence
(Definition A.9 of \cite{FraTheA}). Our approach to this theorem will be in terms
of the Hausdorff metric.

\begin{definition}
Let $(X,d)$ be a compact metric space and let $K_1$ and $K_2$ be nonempty closed subsets of $X$. 
The \textbf{Hausdorff distance} between $K_1$ and $K_2$ is defined to be
\begin{eqnarray*}
d_H(K_1,K_2) 
& = & \max \left\{\sup_{x_1 \in K_1} \inf_{x_2 \in K_2} d(x_1,x_2), \sup_{x_2 \in K_2}
      \inf_{x_1 \in K_1} d(x_1,x_2)\right\}\\
& = & \inf \left\{\varepsilon > 0 |\ K_1 \subseteq N_\varepsilon(K_2) \text{ and }
      K_2 \subseteq N_\varepsilon(K_1)   \right\}
\end{eqnarray*}
where $N_\varepsilon(K) = \bigcup_{y \in K} \{ x \in X |\ d(x,y) \leq \varepsilon \}$. 
\end{definition}

\smallskip\noindent
Note: The Hausdorff distance on the set of all nonempty closed subsets $\mathcal{P}^c(X)$
of a compact metric space $(X,d)$ is a metric, and the two definitions of the Hausdorff
metric given above are equivalent. Moreover, the space $\mathcal{P}^c(X)$ is itself compact
in the topology determined by the Hausdorff metric. (See for instance Section 7.3 
of \cite{MunTop}.)

\smallskip
We would now like to identify a broken flow line with cascades with a closed subset of
some compact metric space. For broken flow lines without cascades this is done by identifying
a broken flow line of a Morse-Bott-Smale function with its image in the compact manifold $M$
(see Section 2 of \cite{HurFlo}). However, a flow line with cascades may have a cascade $x_k$
that ends at a critical point. In this case the parameter $t_k$ records the time spent 
resting at the critical point instead of time spent flowing along the critical submanifold.
Hence, the map that sends a broken flow line with cascades to its image in $M$ is not
injective.  To make this map injective we should keep track of the times $t_k$, in addition
to the image of the broken flow line.

\smallskip
Following \cite{SchMor} we make the following definition.

\begin{definition}
Define the compactification of $\mathbb{R}$ to be $\overline{\mathbb{R}}  =
\mathbb{R} \cup \{\pm \infty\}$ equipped with the structure of a bounded manifold
by the requirement that $\psi:\overline{\mathbb{R}} \rightarrow [-1,1]$ given by
$$
\psi(t) = \frac{t}{\sqrt{1+t^2}}
$$
be a diffeomorphism.
\end{definition}

We also make the following definition regarding the different gradient flows.

\begin{definition}
Let $f:M \rightarrow \mathbb{R}$ be a Morse-Bott function on a Riemannian manifold
$(M,g)$ with critical set $Cr(f) = \coprod_{j=1}^l C_j$, and let 
$f_j:C_j \rightarrow \mathbb{R}$ be a Morse function on the critical submanifold
$C_j$ for $j=1,\ldots ,l$. We define the \textbf{flow} of $\{f,f_1,\ldots ,f_l\}$ on $M$
to be the action $\phi:\mathbb{R} \times M \rightarrow M$ given on a point $x \in M$
for time $t\in \mathbb{R}$ by
$$
\phi_t(x) = \left\{
\begin{array}{ll}
\varphi^f_t(x)     & \text{ if } x \not\in Cr(f) = C_1 \cup \cdots \cup C_l\\
\varphi^{f_j}_t(x) & \text{ if } x \in C_j \text{ for some } j=1,\ldots ,l
\end{array}\right.
$$
where $\varphi^f_t$ denotes the $1$-parameter group of diffeomorphisms generated
by $-\nabla f$ and $\varphi^{f_j}_t$ denotes the $1$-parameter group of diffeomorphisms
generated by $-\nabla f_j$ (with respect to the restriction of $g$ to $C_j$)
for all $j=1,\ldots ,l$.  We extend this action to $\overline{\mathbb{R}}$
by taking limits as $t$ approaches $\pm\infty$.
\end{definition}

\noindent
Note:  The flow of $\{f,f_1,\ldots ,f_l\}$ defines a map $\phi:\overline{\mathbb{R}} 
\times M \rightarrow M$ that is smooth when restricted to $\mathbb{R} \times (M - Cr(f))$
or to $\mathbb{R} \times Cr(f)$.

\smallskip
We now explain how to identify a broken flow line with cascades with an element of 
the compact metric space $\mathcal{P}^c(M) \times\overline{\mathbb{R}}^l$, where $l$
is the number of components of $Cr(f) = \coprod_{j=1}^l C_j$.  Recall that the
space of all nonempty closed subsets of $M$, $\mathcal{P}^c(M)$, is a compact metric
space with respect to the Hausdorff metric. For the metric on $\overline{\mathbb{R}}$ we
will use the totally bounded metric determined by the diffeomorphism 
$\psi:\overline{\mathbb{R}} \rightarrow [-1,1]$. That is, for $x,y \in \mathbb{R}$ we
define
$$
d(x,y) = \left| \frac{x}{\sqrt{1+x^2}} - \frac{y}{\sqrt{1+y^2}}\right| \in [0,2]
$$
and note that $d$ has a unique continuous extension to a metric on $\overline{\mathbb{R}}$.
The space $\mathcal{P}^c(M) \times\overline{\mathbb{R}}^l$ is then a compact
metric space with respect to the product metric.

\smallskip
We will map a broken flow line with cascades $(v_1,\ldots ,v_n)$ to its image 
$\text{Im}(v_1,\ldots ,v_n)$ in $M$ and the time $t_j$ spent flowing along or resting on
each critical submanifold $C_j$ for all $j=1,\ldots ,l$.  This gives a nonempty closed subset
of $M$ and an $l$-tuple in $\overline{\mathbb{R}}^l$, i.e. an element
$\text{Im}(v_1,\ldots ,v_n) \times (t_1,\ldots ,t_l) \in \mathcal{P}^c(M) \times
\overline{\mathbb{R}}^l$.

More explicitly, we define $\text{Im}(v_1,\ldots ,v_n) \subset M$ for a broken flow line
with cascades $(v_1,\ldots ,v_n)$ as follows. Let $\nu \in \{1,\ldots ,n\}$ and suppose
that the unparameterized flow line with $n_\nu$ cascades $v_\nu$ has a parameterization
$$
((x_k^\nu)_{1\leq k \leq n_\nu}, (t_k^\nu)_{1 \leq k \leq n_\nu-1})
$$
where $x_k^\nu \in C^\infty(\mathbb{R},M)$ and $t_k^\nu \in \mathbb{R}_+$.  Then the
image of $v_\nu$ in $M$ is defined to be
$$
\text{Im}(v_\nu) = 
\bigcup_{k=1}^{n_\nu} x_k^\nu(\overline{\mathbb{R}})\ \cup \ \bigcup_{k=1}^{n_\nu -1}
\phi_{[0,t_k^\nu]}(x_k^\nu(\infty)) \subset M
$$
where $\phi_{[0,t_k^\nu]}(x_k^\nu(\infty)) = \bigcup_{0 \leq t \leq t_k^\nu} 
\phi_t(x_k^\nu(\infty))$ and $x_k^\nu(\infty) = \lim_{t \rightarrow \infty} x_k^\nu(t)$.
This definition is clearly independent of the parameterization, and we define 
$\text{Im}(v_1,\ldots ,v_n) \subset M$ to be the union of the images of
$v_\nu$ for all $\nu=1,\ldots ,n$. Note that $\text{Im}(v_1,\ldots ,v_n)$ is the
image of a continuous injective path between two critical points which
is $\overline{\mathbb{R}}$-equivariant with respect to the flow $\phi$ of 
$\{f,f_1,\ldots ,f_l\}$. 

\smallskip
For the other components we map $(v_1,\ldots ,v_n)$ to an $l$-tuple of elements 
$(t_1,\ldots ,t_l) \in \overline{\mathbb{R}}^l$ that records the time spent flowing along or
resting on each critical submanifold.  Explicitly, the $j$th component of this map is 
defined to be:
$$
\begin{array}{rl}
0       & \text{if the image of $(v_1,\ldots ,v_n)$ does not intersect $C_j$}\\
t_j     & \text{if for some $\nu=1,\ldots ,n$ the cascade $v_\nu$ flows along or rests on}\\
        & \text{the critical submanifold $C_j$ for finite time $t_j$}\\
\infty  & \text{otherwise.}
\end{array}
$$
Altogether, this defines an injective map
$$
(v_1,\ldots ,v_n) \mapsto \text{Im}(v_1,\ldots ,v_n) \times (t_1,\ldots ,t_l) \in
\mathcal{P}^c(M) \times\overline{\mathbb{R}}^l.
$$

\begin{definition}
The topology on the space of broken flow lines with cascades is defined by the requirement
that the above injection be a homeomorphism onto its image.
\end{definition}

For $q\in Cr(f_j)$ and $p \in Cr(f_i)$ we will identify the space of broken flow lines
with cascades from $q$ to $p$ with its image under the above injection and
denote this space by $\overline{\mathcal{M}}^c(q,p) \subset \mathcal{P}^c(M) \times
\overline{\mathbb{R}}^l$.

\begin{theorem}\label{compactification}
The space $\overline{\mathcal{M}}^c(q,p)$ is compact, and the injection defined above
restricts to a continuous embedding
$$
\mathcal{M}^c(q,p) \hookrightarrow \overline{\mathcal{M}}^c(q,p) \subset 
\mathcal{P}^c(M) \times \overline{\mathbb{R}}^l.
$$
Hence, every sequence of unparameterized flow lines with cascades from $q$ to $p$ has a 
subsequence that converges to a broken flow line with cascades from $q$ to $p$.
\end{theorem}

\proofstart
Since $\mathcal{P}^c(M) \times\overline{\mathbb{R}}^l$ is compact, any sequence of broken
flow lines $\{(v_1^k,\ldots, v_{n_k}^k)\}$ in $\overline{\mathcal{M}}^c(q,p) \subset
\mathcal{P}^c(M) \times\overline{\mathbb{R}}^l$ must have a subsequence that converges
to some element $C_M \times (t_1,\ldots , t_l) \in \mathcal{P}^c(M) \times
\overline{\mathbb{R}}^l$. We need to show that there exists a subsequence of 
$\{(v_1^k,\ldots, v_{n_k}^k)\}$ (which we still denote by $\{(v_1^k,\ldots, v_{n_k}^k)\}$) 
such that the limit of this subsequence (which we still denote by $C_M \times 
(t_1,\ldots , t_l)$) is in $\overline{\mathcal{M}}^c(q,p) \subset 
\mathcal{P}^c(M) \times\overline{\mathbb{R}}^l$.

We will first show that there exists a subsequence of $\{(v_1^k,\ldots, v_{n_k}^k)\}$ such
that $C_M = Im(v_1,\ldots ,v_n)$ for some broken flow line with cascades $(v_1,\ldots ,v_n)$
from $q$ to $p$. To see this, note that since $Im(v_1^k,\ldots, v_{n_k}^k) \subset M$ is
$\overline{\mathbb{R}}$-equivariant with respect to the flow $\phi$ of $\{f,f_1,\ldots ,f_l\}$
and $\lim_{k \rightarrow \infty} Im(v_1^k,\ldots, v_{n_k}^k) = C_M$ in the Hausdorff metric,
$C_M$ is also $\overline{\mathbb{R}}$-equivariant with respect to the flow $\phi$.
Moreover for every $k$, $Im(v_1^k,\ldots, v_{n_k}^k)$ is the image of a continuous
injective path from $q$ to $p$ with at most one point on each level set $f^{-1}(y)$ 
for every regular value $y$ of $f$ and at most one point on each level set $f_j^{-1}(y)$
for every value $y\in \mathbb{R}$ for all $j=1,\ldots ,l$. Thus, we can pass to a 
subsequence of $\{(v_1^k,\ldots, v_{n_k}^k)\}$ such that the same holds for the limit.
This shows that after passing to an appropriate subsequence we have
$C_M = Im(v_1,\ldots ,v_n)$ for some broken flow line with cascades 
$(v_1,\ldots ,v_n)$ from $q$ to $p$. 

Now let $j \in \{1,\ldots ,l\}$.  For $(t_1,\ldots , t_l)$ there are two cases
to consider: 1) the sequence $\{Im(v_1^k,\ldots, v_{n_k}^k)\}$ does not intersect the
critical submanifold $C_j$ for any $k$ and 2) the sequence $\{Im(v_1^k,\ldots, v_{n_k}^k)\}$
intersects the critical submanifold $C_j$ for all $k$ sufficiently large. Otherwise we can pass
to a subsequence that fits one of these two cases. For the first case, note that the limit 
$C_M$, which is the image of a broken flow line with cascades, can intersect $C_j$ in at 
most one point since $f$ decreases along its gradient flow lines. Thus, for
$\text{Im}(v_1^k,\ldots, v_{n_k}^k) \times (t_1^k,\ldots ,t_l^k) \in 
\mathcal{P}^c(M) \times\overline{\mathbb{R}}^l$ we have $t_j^k=0$ for all $k$,
and $t_j = 0$. For the second case, note that since $\overline{\mathbb{R}}$ is a compact
metric space, we can pass to a subsequence such that $t_j^k \rightarrow t_j$ for some 
$t_j \in \overline{\mathbb{R}}$. By passing to a subsequence for each $j=1,\ldots ,l$ we
obtain an element $(t_1,\ldots , t_l) \in \overline{\mathbb{R}}^l$ such that
$$
\text{Im}(v_1^k,\ldots, v_{n_k}^k) \times (t_1^k,\ldots ,t_l^k) \rightarrow 
Im(v_1,\ldots ,v_n) \times (t_1,\ldots ,t_l) \in \mathcal{P}^c(M) \times\overline{\mathbb{R}}^l
$$
as $k \rightarrow \infty$ and $t_j$ records the time $(v_1,\ldots ,v_n)$ spends flowing
along or resting on each critical submanifold $C_j$ for all $j=1,\ldots ,l$.  
Therefore, every sequence of broken flow lines with cascades from $q$ to $p$ has a
subsequence that converges to a broken flow line with cascades from $q$ to $p$ in
$\overline{\mathcal{M}}^c(q,p) \subset \mathcal{P}^c(M) \times\overline{\mathbb{R}}^l$.

\smallskip
To see that the injection defined above restricts to a continuous embedding
$$
\mathcal{M}^c(q,p) \hookrightarrow \overline{\mathcal{M}}^c(q,p) \subset 
\mathcal{P}^c(M) \times\overline{\mathbb{R}}^l
$$
note that the fibered product and gluing constructions used in the proof
of Theorem \ref{manifold} are compatible with the Hausdorff metric.  That is,
if a sequence of points $v^k$ contained in a smooth chart of $\mathcal{M}^c(q,p)$
converges to a point $v$ in the chart, then
$$
\text{Im}(v^k) \times (t_1^k ,\ldots ,t_l^k) \rightarrow \text{Im}(v) \times
(t_1,\ldots ,t_l)
$$
as $k \rightarrow \infty$.
\proofend

\begin{corollary}\label{finitenumber}
If $\lambda_q - \lambda_p = 1$, then $\mathcal{M}^c(q,p)$ is compact and hence a finite set. 
\end{corollary}

\proofstart
Let $v^k$ be a sequence of unparameterized flow lines with cascades from $q$ to $p$.
By the preceding theorem $v^k$ has a subsequence that converges to a broken flow
line with cascades $(v_1,\ldots ,v_n)$ from $q$ to $p$. Suppose that $v_1$ ends at a
critical point $p'$ with $p' \neq p$. Then Theorem \ref{manifold} implies that
$\lambda_q > \lambda_{p'} > \lambda_p$, which contradicts the assumption that
$\lambda_q - \lambda_p = 1$.  Thus, $p' = p$, $n=1$, and every sequence in $\mathcal{M}^c(q,p)$
has a subsequence that converges to an element of $\mathcal{M}^c(q,p)$. Therefore,
$\mathcal{M}^c(q,p)$ is a compact zero dimensional manifold, i.e. a finite set of points.

\proofend


\smallskip
The preceding corollary allows us to make the following definition under the following
assumptions: 1) $f$ satisfies the Morse-Bott-Smale transversality condition with respect to
the Riemannian metric $g$ on $M$, 2) $f_k:C_k \rightarrow \mathbb{R}$ satisfies the Morse-Smale
transversality condition with respect to the restriction of $g$ to $C_k$ for all 
$k=1,\ldots ,l$, and 3) for all $(i,j)$ and for each pair of critical points $(q,p) \in 
\text{Cr}(f_j) \times \text{Cr}(f_i)$ the unstable and stable manifolds $W^u_{f_j}(q)$ and
$W^s_{f_{i}}(p)$ are transverse to the beginning and endpoint maps. Recall that the total
index of a critical point of $f_j$ was defined in Definition \ref{index} as the Morse index
relative to $f_j$ plus the Morse-Bott index of the critical submanifold $C_j$.

\begin{definition}
Define the $k^{\text{th}}$ chain group $C_k^c(f)$ to be the free abelian group generated 
by the critical points of total index $k$ of the Morse-Smale functions $f_j$ for all
$j=1,\ldots ,l$, and define $n^c(q,p;\mathbb{Z}_2)$ to be the number of flow lines with
cascades between a critical point $q$ of total index $k$ and a critical point $p$ of
total index $k-1$ counted mod 2. Let
$$
C^c_\ast(f) \otimes \mathbb{Z}_2 = \bigoplus_{k=0}^m C_k^c(f) \otimes \mathbb{Z}_2
$$
and define a homomorphism $\partial_k^c:C^c_k(f)\otimes \mathbb{Z}_2 \rightarrow 
C^c_{k-1}(f)\otimes \mathbb{Z}_2$ by
$$
\partial^c_k(q) = \sum_{p \in Cr(f_{k-1})} n^c(q,p;\mathbb{Z}_2)p.
$$
The pair $(C^c_\ast(f) \otimes \mathbb{Z}_2,\partial_\ast^c)$ is called the
\textbf{cascade chain complex} with $\mathbb{Z}_2$ coefficients.
\end{definition}

In the appendix to \cite{FraTheA} there is a continuation theorem that implies that
the cascade chain complex with $\mathbb{Z}_2$ coefficients is, in fact, a chain complex
whose homology is isomorphic to the singular homology $H_\ast(M;\mathbb{Z}_2)$. We will
not prove this here. Instead, we will use the Morse-Smale functions $f_j:C_j \rightarrow
\mathbb{R}$ for $j=1, \ldots ,l$ to define an explicit perturbation of 
$f:M\rightarrow \mathbb{R}$ to a Morse-Smale function $h_\varepsilon:M \rightarrow
\mathbb{R}$ such that for every $k=0,\ldots ,m$
$$
\text{Cr}_k(h_\varepsilon) = \bigcup_{\lambda_j + n = k} \text{Cr}_n(f_j),
$$
where $\lambda_j$ is the Morse-Bott index of the critical submanifold $C_j$.

By proving a correspondence theorem, we will show that for any $q\in \text{Cr}(f_j)$ and
$p\in \text{Cr}(f_i)$ with $\lambda_q - \lambda_p = 1$ there is a one dimensional trivial
cobordism between $\mathcal{M}^c(q,p)$ and $\mathcal{M}_{h_\varepsilon}(q,p)$. This 
cobordism induces an orientation on $\mathcal{M}^c(q,p)$, which allows us to define the
above homomorphism $\partial^c_\ast$ over $\mathbb{Z}$. Moreover, the cobordism shows that
$\partial^c_\ast$ is a boundary operator that agrees with the Morse-Smale-Witten boundary
operator of $h_\varepsilon$ up to sign.


\section{The Correspondence Theorem}\label{correspondencesection}

In this section we define a $1$-parameter family of Morse-Smale functions 
$h_\varepsilon:M \rightarrow \mathbb{R}$ in terms of an explicit perturbation of the 
Morse-Bott-Smale function $f:M\rightarrow \mathbb{R}$. For any $\varepsilon > 0$ the
critical set of $h_\varepsilon$ is given by $\text{Cr}(h_\varepsilon)  = \bigcup_{k=1}^l
\text{Cr}(f_k)$, and the index of a critical point $p \in \text{Cr}(h_\varepsilon)$ agrees
with the total index of $p$.

We prove a correspondence theorem which says that for any $\varepsilon > 0$ sufficiently
small there is a bijection between unparameterized flow lines with cascades and 
unparameterized gradient flow lines of $h_\varepsilon:M \rightarrow \mathbb{R}$ between 
any two critical points $p,q\in \text{Cr}(h_\varepsilon)$ with $\lambda_q - \lambda_p = 1$.
The correspondence theorem allows us to count the number of
unparameterized flow lines with cascades between $q \in \text{Cr}_k(h_\varepsilon)$ and 
$p\in \text{Cr}_{k-1}(h_\varepsilon)$ with sign, which defines an integer $n^c(q,p) \in 
\mathbb{Z}$. 

The integers $n^c(q,p)$ define a homomorphism $\partial^c_k$ analogous to the 
Morse-Smale-Witten boundary operator such that $\partial^c_k = - \partial_k$ (where
$\partial_k$ denotes the Morse-Smale-Witten boundary operator of $h_\varepsilon$).
This shows directly that $\partial^c_{k-1} \circ \partial^c_k = 0$ and the homology
of the cascade chain complex $(C_\ast^c(f), \partial^c_\ast)$ is isomorphic to the homology
of the Morse-Smale-Witten chain complex $(C_\ast(h_\varepsilon), \partial_\ast)$. The
Morse Homology Theorem then implies that the homology of the cascade chain complex with
integer coefficients is isomorphic to the singular homology $H_\ast(M;\mathbb{Z})$.


\subsection{An explicit perturbation}\label{perturbation}

The following perturbation technique, based on \cite{AusMor}, the Morse-Bott Lemma, 
and a folk theorem proved in \cite{AbbLec}, produces an explicit Morse-Smale function
$h_\varepsilon:M \rightarrow \mathbb{R}$ arbitrarily close to a given Morse-Bott-Smale 
function $f:M \rightarrow \mathbb{R}$ such that $h_\varepsilon = f$ outside of a neighborhood
of the critical set $\text{Cr}(f)$. A similar technique was used in \cite{BanDyn} to give a
proof of the Morse-Bott inequalities with somewhat different orientation assumptions than
the classical ``half-space'' method using the Thom Isomorphism Theorem (see \cite{BotMor},
Appendix C of \cite{FarTop}, and Section 2.6 of \cite{NicAnI}).

\medskip
Let $f:M \rightarrow \mathbb{R}$ be a Morse-Bott-Smale function on a finite
dimensional smooth closed Riemannian manifold $(M,g)$.  Let $T_j$ be a small open tubular
neighborhood around each connected component $C_j\subseteq \text{Cr}(f)$ for every 
$j=1,\ldots ,l$ with local coordinates $(u,v,w)$ consistent with those from the 
Morse-Bott Lemma (Lemma \ref{MorseBottLemma}). By ``small" we mean that the following
conditions hold.
\begin{enumerate}
\item Each $T_j$ is contained in the union of the domains of the charts from the Morse-Bott Lemma.
\item For $i\neq j$ we have $T_i \cap T_j = \emptyset$ and $f$ decreases by at least three times
      $\text{max}\{\text{var}(f,T_j)|\ j=1,\ldots, l\}$ along any gradient flow line from $T_i$
      to $T_j$  where $\text{var}(f,T_j) = \sup\{f(x)|\ x \in T_j\} - \inf\{f(x)|\ x \in T_j\}$.
\item If $f(C_i) \neq f(C_j)$, then $\text{var}(f,T_i) + \text{var}(f,T_j) < 
      \left.\left.\frac{1}{3} \right|f(C_i) - f(C_j)\right|$.
\item For every flow line with $n$ cascades between critical points of relative index one
      $((x_k)_{1\leq k \leq n}, (t_k)_{1 \leq k \leq n-1})$, the image of $x_k$ for $k=1,
      \ldots, n$ intersects the closure of exactly two of the tubular neighborhoods 
      $\{T_j\}_{j=1}^l$ (see Definition \ref{flowlinecascade} and Corollary \ref{finitenumber}).
\end{enumerate}
In addition, we will assume that the tubular neighborhoods are small enough 
so that $f:M \rightarrow \mathbb{R}$ still satisfies the Morse-Bott-Smale transversality
condition after modifying the Riemannian metric on the tubular neighborhoods to
make the charts from the Morse-Bott Lemma isometries on $T_j$ with respect to the standard
Euclidean metric on $\mathbb{R}^m$ for all $j=1,\ldots ,l$. From now on we will assume
that the Riemannian metric $g$ has been so modified, i.e. the charts from the Morse-Bott Lemma
are isometries on the tubular neighborhoods with respect to $g$ and the standard 
Euclidean metric on $\mathbb{R}^n$.

\smallskip
Pick positive Morse functions $f_k:C_k \rightarrow \mathbb{R}$ satisfying the 
Morse-Smale transversality condition with respect to the restriction of $g$ to $C_k$ for all 
$k=1,\ldots ,l$ such that for all $i,j = 1, \ldots , l$ and for every pair of critical points
$(q,p) \in \text{Cr}(f_j) \times \text{Cr}(f_i)$ the unstable and stable manifolds
$W^u_{f_j}(q)$ and $W^s_{f_{i}}(p)$ are transverse to the beginning and endpoint maps
(see Lemma \ref{perturb}).  For every $k=1,\ldots ,l$ extend $f_k:C_k\rightarrow \mathbb{R}$
to a function on $T_k$ by making $f_k:T_k \rightarrow \mathbb{R}$ constant in the directions
normal to $C_k$, i.e. $f_k$ is constant in the $v$ and $w$ coordinates coming from the
Morse-Bott Lemma. Let $\tilde{T}_k \subset T_k$ be a smaller open tubular neighborhood of 
$C_k$ with the same coordinates as $T_k$, and let $\rho_k$ be a smooth bump function which is
constant in the $u$ coordinates, equal to $1$ on $\tilde{T}_k$, equal to $0$ outside of 
$T_k$, and strictly decreasing on $T_k - \tilde{T}_k$ with respect to $|v|$ and $|w|$.

Finally, choose $\varepsilon > 0$ small enough so that
$$
\sup_{T_k - \tilde{T}_k} \varepsilon \|\nabla\rho_k f_k\| < \inf_{T_k - \tilde{T}_k}
\|\nabla f\|
$$
for all $k=1,\ldots, l$, and define
$$
h_\varepsilon = f + \varepsilon \left( \sum_{k=1}^l \rho_k f_k \right).
$$
The function $h_\varepsilon:M \rightarrow \mathbb{R}$ is a Morse function close to 
the Morse-Bott-Smale function $f$, and the critical points of $h_\varepsilon$ are
exactly the critical points of the Morse-Smale functions $f_j$ for $j=1,\ldots ,l$.
Moreover, if $q \in C_j$ is a critical point of $f_j:C_j \rightarrow \mathbb{R}$ of 
index $\lambda_q^j$, then $q$ is a critical point of $h_\varepsilon$ of index 
$\lambda_q^{h_\varepsilon} = \lambda_j + \lambda_q^j$, where $\lambda_j$
is the Morse-Bott index of $C_j$.

\begin{lemma}\label{perturbmetric}
There exists an arbitrarily small perturbation of the Riemannian metric $g$
such that $h_{\varepsilon'}:M \rightarrow \mathbb{R}$ is Morse-Smale for all
$0 < \varepsilon' \leq \varepsilon$ with respect to the perturbed metric.
The perturbed metric can be chosen so that it agrees with $g$ on the union of
the tubular neighborhoods $\{T_j\}_{j=1}^l$.
\end{lemma}

\proofstart
Let $\{\varepsilon_i\}_{i=1}^\infty$ be a countable dense subset of $(0,\varepsilon)$.
For every $1 \leq i < \infty$ we can apply Theorem 2.20 of \cite{AbbLec} to conclude that
there is a residual subspace $\mathcal{R}_i$ of the open unit ball $\mathcal{K}_1$ 
in a Banach space $\mathcal{K}$ such that the function $h_{\varepsilon_i}:M \rightarrow
\mathbb{R}$ is Morse-Smale with respect to the Riemannian metric $g + k_i$ for all 
$k_i \in \mathcal{R}_i$. Moreover, we can choose the function $\theta:M \rightarrow [0,\infty)$
in the statement of Theorem 2.20 to be zero on $\bigcup_{j=1}^l T_j$ so that $k_i=0$ on 
$\bigcup_{j=1}^l T_j$ for all $1 \leq i < \infty$.

For any $k \in \bigcap_{i=1}^\infty \mathcal{R}_i$ the Riemannian metric $g+k$ is a metric
that agrees with $g$ on $\bigcup_{j=1}^l T_j$ such that $h_{\varepsilon_i}:M \rightarrow
\mathbb{R}$ is Morse-Smale with respect to $g+k$ for all $1 \leq i < \infty$.  Moreover, 
since $\bigcap_{i=1}^\infty \mathcal{R}_i$ is dense in $\mathcal{K}_1$ we can choose $k \in
\bigcap_{i=1}^\infty \mathcal{R}_i$ arbitrarily close to zero.  This completes the proof of
the lemma since the set of Morse-Smale gradient vector fields is an open and dense subset
of the space of all gradient vector fields on a Riemannian manifold \cite{PalOnM}.

\proofend

Note that we can choose the perturbation of the Riemannian metric small enough so 
that $f:M\rightarrow \mathbb{R}$ still satisfies the Morse-Bott-Smale
transversality condition with respect to the perturbed metric and for all $(i,j)$ and 
for every pair of critical points $(q,p) \in \text{Cr}(f_j) \times \text{Cr}(f_i)$ the
unstable and stable manifolds $W^u_{f_j}(q)$ and $W^s_{f_{i}}(p)$ are still transverse
to the beginning and endpoint maps.

\begin{lemma}\label{constant}
Let $p,q \in \text{Cr}(h_\varepsilon)$ with $\lambda_q - \lambda_p = 1$, and
let  $0 < \varepsilon' \leq \varepsilon$.  If $h_{\varepsilon'}:M \rightarrow \mathbb{R}$
and $h_\varepsilon:M \rightarrow \mathbb{R}$ are Morse-Smale with respect to the same
Riemannian metric, then the number of gradient flow lines of $h_{\varepsilon'}$ from
$q$ to $p$ is equal to the number of gradient flow lines of $h_\varepsilon$ from $q$ to $p$.
\end{lemma}

\proofstart
The lemma will be proved by constructing a one dimensional compact smooth manifold with
boundary $\overline{\mathcal{M}}_{F_{21}}(q,p)$ that is a trivial cobordism between
$\mathcal{M}_{h_{\varepsilon}}(q,p)$ and $\mathcal{M}_{h_{\varepsilon'}}(q,p)$.

Using the notation in Section 6 of \cite{BanMor}, we take $f_1 = h_\varepsilon$, 
$f_2 = h_{\varepsilon'}$, and a smooth homotopy $F_{21}:M \times \mathbb{R} \rightarrow
\mathbb{R}$ that is strictly decreasing in its second component such that for some large
$T \gg 0$ we have
$$
F_{21}(x,t) = \left\{ \begin{array}{llc}
h_\varepsilon(x) - \rho(t) & \text{ if } & t < -T \\
\hat{h}_t(x)   & \text{ if } & -T\leq t \leq T \\
h_{\varepsilon'}(x) - \rho(t) & \text{ if } & t > T,
\end{array}\right.
$$
where $\hat{h}_t(x)$ is an approximation to $\frac{1}{2}(T-t)(h_\varepsilon(x) - \rho(t)) + 
\frac{1}{2}(T+t)(h_{\varepsilon'}(x) - \rho(t))$ that makes $F_{21}$ smooth and
$\rho:\mathbb{R} \rightarrow (-1,1)$ is a smooth strictly increasing function
such that $\lim_{t \rightarrow -\infty} \rho(t) = -1$ and $\lim_{t \rightarrow +\infty}
\rho(t) = 1$. The moduli space of gradient flow lines of $F_{21}:M \times \mathbb{R}
\rightarrow \mathbb{R}$ has a component
$$
\mathcal{M}_{F_{21}}(q,p) = (W^u_{F_{21}}(q) \cap W^s_{F_{21}}(p))/\mathbb{R}
$$
of dimension $1$ (see Lemma 6.2 of \cite{BanMor}) that can be compactified to a smooth
manifold with boundary $\overline{\mathcal{M}}_{F_{21}}(q,p)$ using piecewise
gradient flow lines (see Theorem 6.4 of \cite{BanMor}). 

Moreover, the boundary of the compactified space consists of the fibered products
$$
\partial \overline{\mathcal{M}}_{F_{21}}(q,p) = \overline{\mathcal{M}}_{h_\varepsilon}(q,p)
\times_p \overline{\mathcal{M}}_{F_{21}}(p,p) \coprod \overline{\mathcal{M}}_{F_{21}}(q,q)
\times_q \overline{\mathcal{M}}_{{h_\varepsilon'}}(q,p).
$$
Since $\overline{\mathcal{M}}_{h_\varepsilon}(q,p) \times_p \overline{\mathcal{M}}_{F_{21}}(p,p)
\approx \mathcal{M}_{h_\varepsilon}(q,p)$, $\overline{\mathcal{M}}_{F_{21}}(q,q) \times_q
\overline{\mathcal{M}}_{{h_\varepsilon'}}(q,p) \approx \mathcal{M}_{{h_\varepsilon'}}(q,p)$,
and $F_{21}:M \times \mathbb{R} \rightarrow \mathbb{R}$ is strictly decreasing in
its second component, $\overline{\mathcal{M}}_{F_{21}}(q,p)$ is a one dimensional
trivial cobordism between $\mathcal{M}_{h_\varepsilon}(q,p)$ and 
$\mathcal{M}_{{h_\varepsilon'}}(q,p)$. Thus, $\mathcal{M}_{h_\varepsilon}(q,p)$ and 
$\mathcal{M}_{{h_\varepsilon'}}(q,p)$ have the same number of elements.

\proofend

%

\smallskip\noindent
Remark: The moduli space $\mathcal{M}_{F_{21}}(q,p)$ used in the preceding proof is,
in the language of \cite{SchMor}, a space of $\lambda$-parameterized trajectories between
the trivial regular homotopies $h_\varepsilon$ and $h_{\varepsilon'}$
(see Definition 2.29 of \cite{SchMor}). A general moduli space of $\lambda$-parameterized
trajectories is constructed in Theorem 2 of Section 2.3.2 of \cite{SchMor}, and its compactification is discussed in Section 2.4.4.

\smallskip
In summary, we have a Riemannian metric $g$ on $M$ and a $1$-parameter family of Morse
functions $h_\varepsilon:M \rightarrow \mathbb{R}$ such that the following conditions hold
for all $\varepsilon > 0$ sufficiently small and for all $j=1, \ldots ,l$.

\begin{enumerate}
\item The function $h_0 = f:M \rightarrow \mathbb{R}$ satisfies the Morse-Bott-Smale
      transversality condition with respect to the metric $g$.
\item The functions $h_{\varepsilon}:M \rightarrow \mathbb{R}$ and $f_j:C_j \rightarrow
      \mathbb{R}$ satisfy the Morse-Smale transversality condition with respect to $g$.
\item For all $i,j=1,\ldots ,l$ and for each pair of critical points $(q,p) \in 
      \text{Cr}(f_j) \times \text{Cr}(f_i)$ the unstable and stable manifolds 
      $W^u_{f_j}(q)$ and $W^s_{f_{i}}(p)$ are transverse to the beginning and endpoint maps.
\item The function $h_\varepsilon = f$ outside of the union of the tubular 
      neighborhoods $T_j$.
\item The function $h_\varepsilon = f + \varepsilon f_j$ on the smaller tubular neighborhoods 
      $\tilde{T}_j$. \label{smalltubular}
\item The charts from the Morse-Bott Lemma within the tubular neighborhoods $T_j$ are isometries
      with respect to the metric on $M$ and the standard Euclidean metric on $\mathbb{R}^m$.
      \label{isometries}
\item In the local coordinates $(u,v,w)$ of a tubular neighborhood $T_j$
      we have $f = f(C)- |v|^2 + |w|^2$, $\rho_j$ depends only on the $v$ and $w$
      coordinates, and $f_j$ depends only on the $u$ coordinates. In particular, 
      $\nabla f \perp \nabla f_j$ on $T_j$ by the previous condition. \label{perp}
\item The gradient $\nabla f$ dominates $\varepsilon \nabla \rho_j f_j$ on 
      $T_j - \tilde{T}_j$. \label{dominate}
\item For $q,p \in \text{Cr}(h_\varepsilon)$ with $\lambda_q - \lambda_p = 1$,
      the number of gradient flow lines of $h_\varepsilon$ from $q$ to $p$ is
      independent of $\varepsilon > 0$.
\end{enumerate}

\begin{lemma}\label{degenerating}
Let $\varepsilon > 0$ be small enough so that the above conditions hold, and let
$\{\varepsilon_\nu\}_{\nu=1}^\infty$ be a decreasing sequence such that $0 < \varepsilon_\nu
\leq \varepsilon$ for all $\nu$ and $\lim_{\nu \rightarrow \infty} \varepsilon_\nu = 0$. Let 
$q,p \in \text{Cr}(h_\varepsilon)$, and suppose that $\gamma_{\varepsilon_\nu} \in
\mathcal{M}_{h_{\varepsilon_\nu}}(q,p)$ for all $\nu$. Then there exists a broken flow line
with cascades $\gamma \in \overline{\mathcal{M}}^c(q,p)$ and a subsequence of
$\{Im(\gamma_{\varepsilon_\nu})\}_{\nu=1}^\infty$ that converges to $Im(\gamma)$ in the
Hausdorff topology.
\end{lemma}

\proofstart
Let $q \in C_j$, $p \in C_i$, and $\gamma_{\varepsilon_\nu} \in \mathcal{M}_{h_{\varepsilon_\nu}}
(q,p)$ where $\lim_{\nu \rightarrow \infty} \varepsilon_\nu = 0$. Recall that outside of the
open tubular neighborhoods $\{T_k\}_{k=1}^l$ we have $h_{\varepsilon_\nu} = f$, and inside
$T_k$ we have
$$
h_{\varepsilon_\nu} = f + \varepsilon_\nu \rho_k f_k
$$
where $\nabla f \perp \nabla f_k$, $0 \leq \rho_k \leq 1$, and $f_k > 0$. Moreover,
$\nabla h_{\varepsilon_\nu} = \nabla f + \varepsilon_\nu \nabla f_k$ on the smaller open tubular
neighborhood $\tilde{T_k} \subset T_k$, and $\nabla f$ dominates $\varepsilon_\nu \nabla \rho_k f_k$
on $T_k - \tilde{T}_k$. By passing to a subsequence of $\{\gamma_{\varepsilon_\nu}\}_{\nu=1}^\infty$
we may assume that there exists a set of distinct integers $\{j_1, j_2, \ldots, j_{n-1}\} \subseteq 
\{1, 2, \ldots, l\}$ such that for all $\nu$ we have $\text{Im}(\gamma_\nu) \cap T_{j_k} \neq 
\emptyset$ for all $k=1, \ldots, n-1$ and $\text{Im}(\gamma_\nu) \cap T_k = \emptyset$ if $k \in 
\{1, 2, \ldots, l\} - \{i, j_1, j_2, \ldots, j_{n-1}, j\}$.
 
Since $\mathcal{P}^c(M)$ is compact in the Hausdorff topology, there exists a
subsequence of $\{\gamma_{\varepsilon_\nu}\}_{\nu=1}^\infty$, which we still denote by
$\{\gamma_{\varepsilon_\nu}\}_{\nu=1}^\infty$, such that the compact sets
$$
C_{\varepsilon_\nu} = \text{Im}(\gamma_{\varepsilon_\nu}) - \left( T_i \cup 
\bigcup_{k=1}^{n-1} T_{j_k} \cup T_j \right)
$$
converge to some compact set $C \in \mathcal{P}^c(M)$ as $\nu \rightarrow \infty$.
The interior of each $C_{\varepsilon_\nu}$ is locally invariant under the flow of $-\nabla f$,
and hence the interior of the limit $C$ is also locally invariant with respect to the flow of
$-\nabla f$.  Moreover, for every regular value $y$ of $f$ the level set $f^{-1}(y)$
contains at most one element of $C_{\varepsilon_\nu}$ for each $\nu$, and hence we can
pass to a subsequence of $\{\gamma_{\varepsilon_\nu}\}_{\nu=1}^\infty$ such that the same
holds for $C$. Therefore, there exists a subsequence of
$\{\gamma_{\varepsilon_\nu}\}_{\nu=1}^\infty$, which we still denote by 
$\{\gamma_{\varepsilon_\nu}\}_{\nu=1}^\infty$, and gradient flow lines $x_1, \ldots, x_n$
of $-\nabla f$ (not necessarily distinct) such that 
$$
\text{Im}(\gamma_{\varepsilon_\nu}) - \left( T_i \cup \bigcup_{k=1}^{n-1} T_{j_k} 
\cup T_j \right) \rightarrow
\bigcup_{k=1}^n\text{Im}(x_k) - \left( T_i \cup \bigcup_{k=1}^{n-1} T_{j_k} \cup T_j \right)
$$
in the Hausdorff topology as $\nu \rightarrow \infty$.  Moreover, since 
$\nabla h_{\varepsilon_\nu} = \nabla f + \varepsilon_\nu \nabla \rho_k f_k$ and there is a 
positive lower bound for $\|\nabla f\|$ on $T_k - \tilde{T}_k$ for all $k=1,\ldots , l$
we have
$$
\text{Im}(\gamma_{\varepsilon_\nu}) - \left( \tilde{T}_i \cup \bigcup_{k=1}^{n-1}
\tilde{T}_{j_k} \cup \tilde{T}_j \right) \rightarrow
\bigcup_{k=1}^n\text{Im}(x_k) - \left( \tilde{T}_i \cup \bigcup_{k=1}^{n-1} \tilde{T}_{j_k} 
\cup \tilde{T}_j \right)
$$
in the Hausdorff topology as $\nu \rightarrow \infty$. We will order the gradient flow lines
$x_1, \ldots, x_n$ as in Definition \ref{flowlinecascade}, i.e. $x_k(t)$ flows into 
$T_{j_k}$ as $t$ increases for all $k=1, \ldots, n-1$.

On the tubular neighborhood $\tilde{T}_j$ we have $\nabla h_{\varepsilon_\nu} = \nabla f +
\varepsilon_\nu \nabla f_j$ where $\nabla f \perp \nabla f_j$, and hence there is a subsequence of 
$\{\gamma_{\varepsilon_\nu}\}_{\nu=1}^\infty$ such that $\text{Im}(\gamma_{\varepsilon_\nu})
\cap \tilde{T}_j$ converges to a curve consisting of the union of $\text{Im}(x_1) \cap
\tilde{T}_j$ and a (possibly broken) gradient flow line of $f_j$ from 
$q$ to $\lim_{t \rightarrow -\infty} x_1(t)$. Similar statements apply to the tubular neighborhood
$\tilde{T}_i$.

For each tubular neighborhood $\tilde{T}_{j_1}, \ldots , \tilde{T}_{j_{n-1}}$ there are
two cases to consider: 1) there exists a neighborhood $U \subseteq \tilde{T}_{j_k}$ of 
$C_{j_k}$ such that $\text{Im}(\gamma_{\varepsilon_\nu}) \cap U = \emptyset$ for
all $\nu$ or 2) for every neighborhood $U \subseteq \tilde{T}_{j_k}$ of $C_{j_k}$ we have
$\text{Im}(\gamma_{\varepsilon_\nu}) \cap U \neq \emptyset$ for all $\nu$ sufficiently large.
Otherwise we can pass to a subsequence of $\{\gamma_{\varepsilon_\nu}\}_{\nu=1}^\infty$ such
that one of these cases applies. In the first case, there is a positive lower bound for
$\|\nabla f\|$ on $\text{Im}(\gamma_{\varepsilon_\nu}) \cap \tilde{T}_{j_k}$ independent of $\nu$,
and hence $\nabla h_{\varepsilon_\nu}$ converges to $\nabla f$ on $\text{Im}(\gamma_{\varepsilon_\nu})
\cap \tilde{T}_{j_k}$ as $\nu \rightarrow \infty$. Thus, $x_k(t)$ and $x_{k+1}(t)$ are the
same gradient flow line of $f$, and $\text{Im}(\gamma_{\varepsilon_\nu}) \cap \tilde{T}_{j_k}$
converges to $\text{Im}(x_k) \cap \tilde{T}_{j_k}$ as $\nu \rightarrow \infty$.

In the second case, $\lim_{t \rightarrow \infty} x_k(t) \in C_{j_k}$ since 
$\text{Im}(\gamma_{\varepsilon_\nu}) \cap f^{-1}(y)$ converges to $\text{Im}(x_k) \cap f^{-1}(y)$ 
for any $y > f(C_{j_k})$ with $\text{Im}(\gamma_{\varepsilon_\nu}) \cap f^{-1}(y) \in \tilde{T}_{j_k}$.
Similarly, $\lim_{t \rightarrow -\infty} x_{k+1}(t) \in C_{j_k}$. Moreover, 
$\text{Im}(\gamma_{\varepsilon_\nu}) \cap \tilde{T}_{j_k}$ converges to the union of
$\text{Im}(x_k) \cap \tilde{T}_{j_k}$, $\text{Im}(x_{k+1}) \cap \tilde{T}_{j_k}$ and a curve in 
$C_{j_k}$ from $\lim_{t \rightarrow \infty} x_k(t)$ to $\lim_{t \rightarrow -\infty}
x_{k+1}(t)$. Since $\nabla h_{\varepsilon_\nu} = \nabla f + \varepsilon_\nu \nabla f_{j_k}$ in
$\tilde{T}_{j_k}$ where $\nabla f \perp \nabla f_{j_k}$, the curve in $C_{j_k}$ must be a
subset of the image of a (possibly broken) gradient flow line of $f_{j_k}$.
Therefore, there exists a subsequence of $\{\gamma_{\varepsilon_\nu}\}_{\nu=1}^\infty$ and a
broken flow line with cascades $\gamma \in \overline{\mathcal{M}}^c(q,p)$ such that 
$\{Im(\gamma_{\varepsilon_\nu})\}_{\nu=1}^\infty$ converges to $Im(\gamma)$ in the Hausdorff topology.

\proofend


\subsection{Correspondence theorem}

Throughout this subsection we will assume that the function
$$
h_\varepsilon = f + \varepsilon \left( \sum_{k=1}^l \rho_k f_k \right)
$$
and the Riemannian metric $g$ on $M$ satisfy all the conditions listed
above.  The main goal of this subsection is to prove the following.

\begin{theorem}[Correspondence of Moduli Spaces]\label{correspondence}
Let $p,q \in \text{Cr}(h_\varepsilon)$ with $\lambda_q - \lambda_p = 1$.
For any sufficiently small $\varepsilon>0$ there is a bijection between unparameterized
cascades and unparameterized gradient flow lines of the Morse-Smale function
$h_\varepsilon:M \rightarrow \mathbb{R}$ between $q$ and $p$,
$$
\mathcal{M}^c(q,p) \leftrightarrow \mathcal{M}_{h_\varepsilon}(q,p).
$$
\end{theorem}

\medskip
We will prove this theorem using results from geometric singular perturbation theory
\cite{JonGeo}.  In particular, we will use the Exchange Lemma for fast-slow systems
\cite{JonGen} \cite{SchExcI} \cite{SchExcII}. Roughly speaking, the Exchange Lemma says 
that a manifold $M_0$ that is transverse to the stable manifold of a normally hyperbolic
locally invariant submanifold $C$ will have subsets that flow forward in time under the
full fast-slow system to be near subsets of the unstable manifold of $C$. The Exchange
Lemma can be viewed as a generalization of the $\lambda$-Lemma, which applies to hyperbolic
fixed points (see for instance Theorem 6.17 and Corollary 6.20 of \cite{BanLec}).

In our setup, we have tubular neighborhoods $T_j$ of the critical submanifolds $C_j$ 
for all $j=1,\ldots, l$ and local coordinate charts on $T_j$ that are isometries with
respect to the standard Euclidean metric on $\mathbb{R}^m$. We also have smaller
tubular neighborhoods $\tilde{T}_j \subset T_j$ such that within the smaller tubular
neighborhoods the negative gradient flow of $h_\varepsilon:T \rightarrow \mathbb{R}$
constitutes a fast-slow system because $\nabla h_\varepsilon = \nabla f + \varepsilon
\nabla f_j$ and $\nabla f \perp \nabla f_j$. Moreover, we have coordinates $(u,v,w)$ 
where the function $f|_{\tilde{T_j}}$ depends only on the $(v,w)$ coordinates, 
which are the fast variables, and the function $f_j|_{\tilde{T_j}}$ depends only on 
the $u$ variables, which are the slow variables.

\medskip\noindent
\textbf{Proof of Theorem \ref{correspondence}:} Let $q \in \text{Cr}(f_j)$ and $p \in
\text{Cr}(f_i)$. An unparameterized cascade $\gamma \in \mathcal{M}^c(q,p)$ can be represented by a
flow line with $n$ cascades from $q$ to $p$: $((x_k)_{1 \leq k \leq n}, (t_k)_{1 \leq k \leq n-1})$,
where $t_k$ is the time spent flowing along (or resting on) the intermediate critical submanifold
$C_{j_k}$. For $1 \leq k \leq n-1$, let $y_k:\mathbb{R} \rightarrow C_{j_k}$ be the parameterized
gradient flow line of $f_{j_k}:C_{j_k} \rightarrow \mathbb{R}$ satisfying $y_k(0) = 
\lim_{t \rightarrow \infty} x_k(t)$ and $y_k(t_k) = \lim_{t \rightarrow -\infty} x_{k+1}(t)$
(as in Definition \ref{flowlinecascade}). Assume that $y_k(0) \neq y_k(t_k)$ for any
$1 \leq k \leq n-1$.  This last condition is required in order to apply the Exchange Lemma, and
it holds whenever $\lambda_q - \lambda_p = 1$. To see this, note that if $y_k(0) = y_k(t_k)$
then there is a piecewise gradient flow line of $f$ from the beginning of $x_k$ to the end of
$x_{k+1}$. Hence, there is a 1-parameter family of gradient flow lines of $f$ from the
beginning of $x_k$ to the end of $x_{k+1}$ by the gluing theorem for Morse-Bott moduli spaces
(see the proof of Theorem \ref{manifold}). Each of these gradient flow lines determines a
unique flow line with cascade from $q$ to $p$, and hence $\text{dim }\mathcal{M}^c(q,p) \geq 1$.

For every $1 \leq k \leq n-1$, let $S_k \subset C_{j_k}$ be a tubular neighborhood
of the image $y_k([0,t_k])$ that is diffeomorphic to some contractible open subset
$U_k \subset \mathbb{R}^{\text{dim } C_{j_k}}$. The tubular neighborhood $S_k$ exists
because $y_k([0,t_k])$ is contractible and hence has a trivial normal bundle in $C_{j_k}$.
Similarly, the normal bundle of $S_k \subset M$ is trivial, and hence $S_k$ has a contractible
tubular neighborhood in $\tilde{T}_{j_k}$. This establishes Fenichel coordinates $(u,v,w)$
near $S_k$. (See Proposition 1 and Section 6 of \cite{JonGen}, but note that we do not need
$S_k$ to vary with $\varepsilon$.)

Let $B^k_{\Delta,U_k}$ be a small ``box'' in the phase space $\mathbb{R}^m$ with respect
to the Fenichel coordinates near $S_k$, e.g. 
$$
B^k_{\Delta,U_k} = \{(u,v,w) \in \mathbb{R}^m|\ |v| < \Delta,\ |w| < \Delta,\ u \in U_k\} 
$$
for some small $\Delta > 0$, and let $B_k$ denote the image of $B^k_{\Delta,U_k}$
in $M$. We will show that for $\Delta > 0$ and $\varepsilon > 0$ sufficiently small
there exist submanifolds $M_k \subset W^u_{h_\varepsilon}(q)$ that satisfy the following
conditions for every $1 \leq k \leq n-1$.
\begin{itemize}
  \item[(\text{D}1)]   $\lambda_{j_k} \leq \text{dim } M_k \leq \lambda_{j_k} + 
                       \text{dim }C_{j_k} - 1$
  \item[$(\text{T}1)$] There exists a point $q_k \in M_k \cap \overline{B}_k$ such that $M_k
                       \pitchfork_{q_k} W^s_f(S_k)$.
  \item[$(\text{T}2)$] The omega limit set $J_k = \omega(M_k \cap W^s_f(S_k) \cap V_k) \subset
                       S_k$ with respect to the flow of $-\nabla f$ is a manifold of dimension
                       $\text{dim } M_k - \lambda_{j_k}$, where $V_k$ is a small enough
                       open neighborhood of $q_k$ to ensure that $M_k \cap W^s_f(S_k) \cap
                       V_k$ is a manifold, and $\nabla f_{j_k}$ is not tangent to $J_k$.
  \item[$(\text{T}3)$] The tangent space to $M_k$ at $q_k$ intersects the tangent space
                       of $W^s_f(\omega(q_k))$ in a zero dimensional space.
  \item[$(\text{I}1)$] If $\text{Im}(\gamma_\varepsilon) \cap M_k \neq \emptyset$ for some 
                       $\gamma_\varepsilon \in \mathcal{M}_{h_\varepsilon}(q,p)$,
                       then $\text{Im}(\gamma'_\varepsilon) \cap M_k \neq \emptyset$ 
                       for every $\gamma'_\varepsilon \in \mathcal{M}_{h_\varepsilon}(q,p)$
                       with $d_H(\text{Im}(\gamma'_\varepsilon), \text{Im}(\gamma)) \leq
                       d_H(\text{Im}(\gamma_\varepsilon), \text{Im}(\gamma))$.
\end{itemize}


The manifold $M_1$ exists as long as $\varepsilon > 0$ is small enough so that the conditions
listed in the previous subsection hold. That is, the conditions in the previous subsection
imply that $\lim_{t \rightarrow -\infty} x_1(t) \in W^u_{h_\varepsilon}(q)$ and 
$W^u_f(\lim_{t \rightarrow -\infty} x_1(t)) \pitchfork W^s_f(S_1)$. Thus, we can find a small 
open neighborhood in $W^u_{h_\varepsilon}(q)$ around the point $r_1$ where the image of $x_1$ 
intersects the boundary of $T_j$ with a cross section that intersects $W^s_f(S_1)$ transversally.
This cross section flows forward under the flow of $-\nabla h_\varepsilon$ to a submanifold
$\tilde{M}_1$ of dimension $\lambda_q - 1$ that intersects $\overline{B}_1 \cap W^s_f(S_1)$
at some point $q_1$. The Morse-Bott-Smale transversality condition implies that $\lambda_{j_1}
< \lambda_j$ (see Lemma 3.6 of \cite{BanMor}), and hence $\lambda_{j_1} \leq \lambda_j +
\lambda_q^j - 1 = \lambda_q - 1 = \text{dim }\tilde{M}_1$. If $\text{dim }\tilde{M}_1
\leq \lambda_{j_1} + \text{dim } C_{j_1} - 1$, then we can take $M_1 = \tilde{M}_1$. 
Otherwise, we can find a small open ball $M_1 \subset \tilde{M}_1$ of dimension 
$\lambda_{j_1} + \text{dim }C_{j_1} - 1$ that satisfies the above conditions. Thus,
$M_1$ exists and $\text{dim }M_1 = \text{min} \{\lambda_q - 1, \text{dim }C_{j_1} + 
\lambda_{j_1} - 1\}$.

\smallskip
We will see by induction using the Exchange Lemma that for $\Delta > 0$ and $\varepsilon > 0$
sufficiently small $M_k \subseteq W^u_{h_\varepsilon}(q)$ exists for $k=2,\ldots ,n-1$.
For this purpose, assume that $\Delta > 0$ and $\varepsilon > 0$ are small enough
so that the conditions listed in the previous subsection hold, the Exchange Lemma applies
around $S_k$ for all $k=1, \ldots , n-1$, and $M_1$ exists.
%
Assume that for some $k$ there exists a submanifold $M_k \subseteq W^u_{h_\varepsilon}(q)$ 
that satisfies the above conditions, and let $M_k^\ast$ and $J_k^\ast$ denote the manifolds
obtained by flowing $M_k$ and $J_k$ forward in time with respect to $-\nabla h_\varepsilon$
on the time interval $[0,\infty)$. The dimension of $M_k^\ast$ is $\text{dim }M_k + 1$,
and $\text{dim }J_k^\ast = \text{dim } M_k^\ast - \lambda_{j_k}$.

Let $x^\varepsilon_{k+1}(t)$ be the gradient flow line of $h_\varepsilon$ through the point
$r_{k+1}$ where the image of $x_{k+1}(t)$ intersects the boundary of $T_{j_k}$. We have
$\lim_{t \rightarrow -\infty} x^\varepsilon_{k+1}(t) = \lim_{t \rightarrow -\infty} y_k(t)$.
Hence, as long as $\varepsilon > 0$ is sufficiently small, the point where
$x_{k+1}^\varepsilon(t)$ exits the box $B_k$ will be in $W^u_f(J_k^\ast)$. Choose a small
open disk $D_k$ in $W^u_f(J_k^\ast)$ of dimension $\text{dim }M_k^\ast$ around this point.
The Exchange Lemma implies that by decreasing $\varepsilon > 0$ we can find an open disk 
$\tilde{D}_k$ in $M_k^\ast$ as close as we like to $D_k$. (See for instance Theorem 6.5 of 
\cite{JonGen}, Lemma 6 of \cite{JonGeo}, or Theorem 2.3 of \cite{SchExcII}.) In this context
``close'' can be in the sense of Definition 6.13 of \cite{BanLec} or ``close'' in the sense
that $\tilde{D}_k$ can be expressed as the graph of a vector valued function over $D_k$ that
goes to zero exponentially along with its derivatives up to finite order as
$\varepsilon \rightarrow 0$ \cite{SchExcII}.

\begin{figure}
\includegraphics{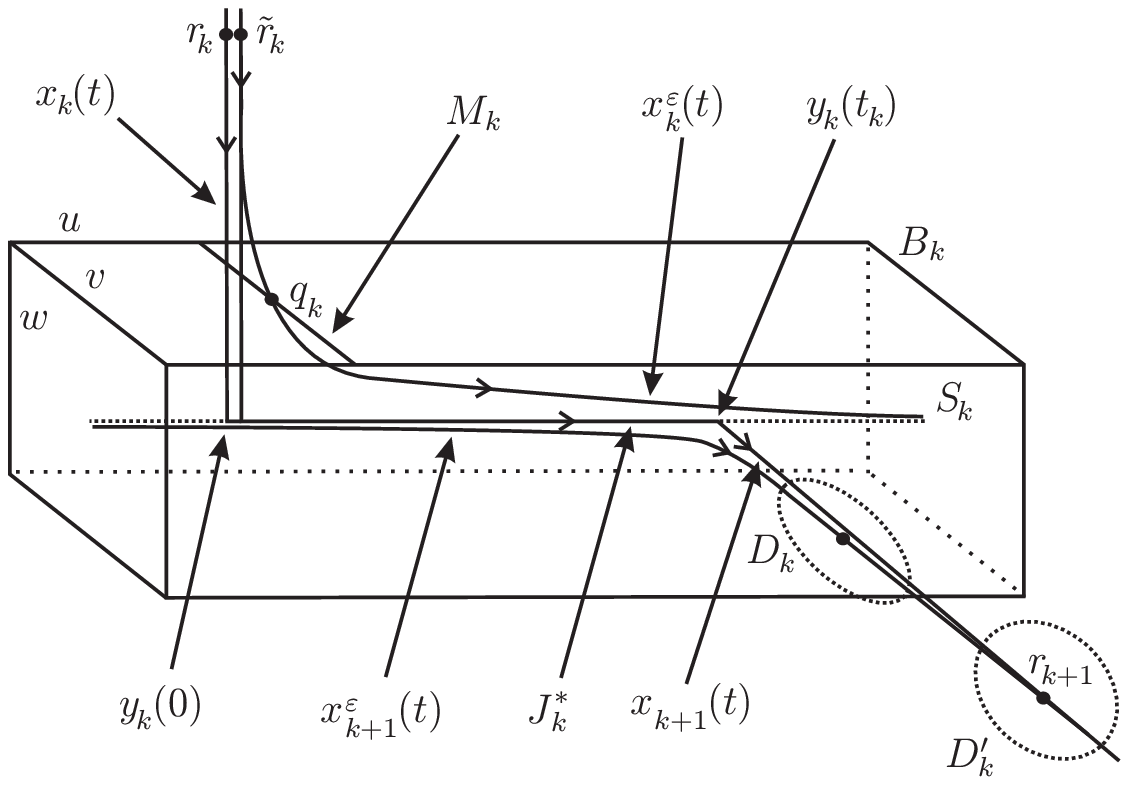}
\end{figure}

The open disk $D_k$ flows forward in finite time under the flow of $-\nabla h_\varepsilon$
to a neighborhood $D_k'$ of $r_{k+1}$, and the open disk $\tilde{D}_k$ flows forward under
the same flow to an open set $\tilde{D}_k'$ in $M_k^\ast$ close to $D_k'$. In fact,
inside $T_{j_k}$ the disks $\tilde{D}_k$ and $D_k$ get closer under the forward time flow
of $-\nabla h_\varepsilon$. The Morse-Bott-Smale transversality condition implies that
$D_k' \pitchfork W^s_f(S_{k+1})$, and hence $\tilde{D}_k' \pitchfork W^s_f(S_{k+1})$ if
$\tilde{D}_k'$ is close enough to $D_k'$ since the collection of maps transverse to a given
submanifold is locally stable (see for instance Theorem 5.16 of \cite{BanLec} or Theorem
3.2.1 of \cite{HirDif}). Thus, we can decrease $\varepsilon > 0$, if necessary, to obtain
an open set $\tilde{D}_k' \subset M_k^\ast$ such that $\tilde{D}_k' \pitchfork W^s_f(S_{k+1})$.
Moreover, $r_{k+1} \in D_k' \cap W^s_f(S_{k+1}) \neq \emptyset$, and hence there exists a
point $\tilde{r}_{k+1} \in \tilde{D}_k' \cap W^s_f(S_{k+1})$ such that $\tilde{D}_k'
\pitchfork_{\tilde{r}_{k+1}} W^s_f(S_{k+1})$.

For $\varepsilon > 0$ sufficiently small, the point $\tilde{r}_{k+1} \in W^s_f(S_{k+1})$
flows forward in time under the flow of $-\nabla h_\varepsilon$ to a point $q_{k+1} \in
\partial \overline{B}_{k+1}$ since the tubular neighborhoods $\{T_j\}_{j=1}^l$ were chosen
small enough so that the image of $x_{k+1}$ does not intersect the closure of any of the
tubular neighborhoods other than $\overline{T}_{j_{k}}$ and $\overline{T}_{j_{k+1}}$.
Moreover, $\tilde{D}_k' \subseteq M_k^\ast$ flows forward in time under the flow of 
$-\nabla h_\varepsilon$ to a submanifold of $W^u_{h_\varepsilon}(q)$ that is transverse to
$W^s_f(S_{k+1})$ at $q_{k+1}$. Thus, we can find a manifold $M_{k+1} \subset M_k^\ast \subset
W^u_{h_\varepsilon}(q)$ of dimension $\text{min}\{\text{dim }M_k, \text{dim }C_{j_{k+1}} +
\lambda_{j_{k+1}} - 1\}$ that satisfies the above conditions. This completes the induction.
Note that if we have to decrease $\varepsilon > 0$ during the induction, then we also have
to modify $M_k \subset W^u_{h_\varepsilon}(q)$. However, $\varepsilon > 0$ will only
need to be decreased a finite number of times.  Hence, we can find a sufficiently small
$\varepsilon > 0$ so that $M_k \subset W^u_{h_\varepsilon}(q)$ exists for all
$k=1, \ldots, n-1$.

\smallskip
To summarize, we have shown that for $\Delta > 0$ and $\varepsilon > 0$ sufficiently
small there exist submanifolds $M_k \subseteq  W^u_{h_\varepsilon}(q)$ and points $q_k$
such that $M_k \pitchfork_{q_k} W^s_f(S_k)$ for all $k=1, \ldots , n-1$. Moreover, 
the point $q_k$ is the image under the forward time flow of $-\nabla h_\varepsilon$ 
of a point $\tilde{r}_k \in W^s_f(S_k) \cap W^u_{h_\varepsilon}(q)$ close to the point
$r_k$ where the image of $x_k(t)$ intersects the boundary of $T_{j_{k-1}}$, 
$$
M_{n-1}^\ast \subset M_{n-2}^\ast \subset \cdots \subset M_2^\ast \subset M_1^\ast
\subset W^u_{h_\varepsilon}(q), 
$$
and every gradient flow line in $\mathcal{M}_{h_\varepsilon}(q,p)$ whose image
is sufficiently close to the image of the cascade $\gamma \in \mathcal{M}^c(q,p)$
intersects $M_{n-1}^\ast$ (and hence $M_k$ for all $k=1, \ldots, n-1$). We can now
repeat the above argument involving the Exchange Lemma for $M_{n-1}$ to see that for
$\Delta > 0$ and $\varepsilon > 0$ sufficiently small we can find an open neighborhood 
$\tilde{D}_{n-1}' \subset M_{n-1}^\ast$ as close as we like to a small open neighborhood
$D_{n-1}'\subset W^u_f(J^\ast_{n-1})$ around the point $r_n$ where the image of
$x_n(t)$ intersects the boundary of $T_{j_{n-1}}$.

Now recall the assumption that
$$
\mathcal{M}_n^c(C_j,C_{j_1}, \ldots , C_{j_{n-1}}, C_i) \stackrel{(\partial_-,\partial_+)}
{\longrightarrow} C_j \times C_i
$$
is transverse and stratum transverse to $W^u_{f_j}(q) \times W^s_{f_i}(p)$ 
(Definition \ref{beginningendtransverse}). This implies that
$$
\mathcal{M}_f(J_{n-1}^\ast,C_i) \stackrel{\partial_+}{\longrightarrow} C_i
$$
is transverse to $W^s_{f_i}(p)$ at $\lim_{t \rightarrow \infty}x_n(t)\in W^s_{f_i}(p)$,
since the endpoint map $\partial_+: \mathcal{M}_n^c(W^u_{f_j}(q),C_{j_1}, \ldots ,
C_{j_{n-1}}, C_i) \rightarrow C_i$ factors through $\partial_+:\mathcal{M}_f
(J_{n-1}^\ast,C_i)\rightarrow C_i$ and is transverse to $W^s_{f_i}(p)$ at $\lim_{t
\rightarrow \infty}x_n(t)$.
%
%
Therefore, $D'_{n-1} \pitchfork_{r_n} W^s_f(W^s_{f_i}(p))$ as long as $D'_{n-1}$ is
sufficiently small. Thus if $\varepsilon > 0$ is sufficiently small, there exists a point
$\tilde{r}_n \in\tilde{D}'_{n-1} \cap W^s_f(W^s_{f_i}(p))$ close to $r_n$ such that
$\tilde{D}'_{n-1} \pitchfork_{\tilde{r}_n} W^s_f(W^s_{f_i}(p))$. The unparameterized
gradient flow line of $h_\varepsilon$ that passes through $\tilde{r}_n$ is an element
$\gamma_{\tilde{r}_n} \in \mathcal{M}_{h_\varepsilon}(q,p)$ whose image is close to the
image of the cascade in $\mathcal{M}^c(q,p)$ represented by $((x_k)_{1 \leq k \leq n},
(t_k)_{1 \leq k \leq n-1})$.

Moreover, if $\lambda_q - \lambda_p = 1$ then we can choose the $D'_k$ small enough
so that $\gamma \in \mathcal{M}^c(q,p)$ is the unique element whose image intersects
$D'_k$ for all $k=1, \ldots, n-1$. Then if $\tilde{D}'_{n-1}$ is sufficiently close to
$D_{n-1}'$, the gradient flow line of $h_\varepsilon$ through $\tilde{r}_n$ will be the
unique element of $\mathcal{M}_{h_\varepsilon}(q,p)$ whose image intersects 
$\tilde{D}_{n-1}' \subset M_{n-1}^\ast$. Thus for $\lambda_q - \lambda_p = 1$ and 
$\varepsilon > 0$ sufficiently small we have defined an injective map
$$
\mathcal{M}^c(q,p) \rightarrow \mathcal{M}_{h_\varepsilon}(q,p)
$$
that sends a cascade $\gamma \in \mathcal{M}^c(q,p)$ to a gradient flow line 
$\gamma_\varepsilon \in \mathcal{M}_{h_\varepsilon}(q,p)$ such that $\text{Im}(\gamma)$ is
close to $\text{Im}(\gamma_\varepsilon)$ in the Hausdorff topology. To see that this
map is surjective, first recall that Lemma \ref{constant} says that if $\varepsilon > 0$
is sufficiently small, then the (finite) number of elements in 
$\mathcal{M}_{h_\varepsilon}(q,p)$ does not depend on $\varepsilon > 0$. So, if the above
map were not surjective, we could pick a decreasing sequence 
$\{\varepsilon_\nu\}_{\nu=1}^\infty$ with $\lim_{\nu \rightarrow \infty} \varepsilon_\nu = 0$
and a sequence of elements $\gamma_{\varepsilon_\nu} \in \mathcal{M}_{h_{\varepsilon_\nu}}(q,p)$ 
such that $\gamma_{\varepsilon_\nu}$ is not in the image of the map 
$$
\mathcal{M}^c(q,p) \rightarrow \mathcal{M}_{h_\varepsilon}(q,p) \leftrightarrow
\mathcal{M}_{h_{\varepsilon_\nu}}(q,p)
$$
for all $\nu$. Lemma \ref{degenerating} would then imply that there exists a subsequence of 
$\{Im(\gamma_{\varepsilon_\nu})\}_{\nu=1}^\infty$ (which we still denote by 
$\{Im(\gamma_{\varepsilon_\nu})\}_{\nu=1}^\infty$) that converges to the image of some
element $\gamma \in \overline{\mathcal{M}}^c(q,p) = \mathcal{M}^c(q,p)$ in the Hausdorff
topology. But if we were to apply the above construction to $\gamma$, then for $\nu$ sufficiently
large we would get an element $\gamma_{\tilde{r}_n} \in \mathcal{M}_{h_{\varepsilon_\nu}}(q,p)$ that 
intersects an open neighborhood $\tilde{D}_{n-1}' \subset M^\ast_{n-1}$ near $W^u_f(C_{j_{n-1}})$.
Since the sequence $\{Im(\gamma_{\varepsilon_\nu})\}_{\nu=1}^\infty$ is converging to
$\text{Im}(\gamma)$  we must have $\text{Im}(\gamma_{\varepsilon_\nu}) \cap \tilde{D}'_{n-1}
\neq \emptyset$ for $\nu$ sufficiently large by condition $(\text{I}1)$, and since 
$\gamma_{\tilde{r}_n}$ is the unique gradient flow line in $\mathcal{M}_{h_{\varepsilon_\nu}}(q,p)$
whose image intersects $\tilde{D}'_{n-1}$, we see that $\gamma_{\varepsilon_\nu} = 
\gamma_{\tilde{r}_n}$ is in the image of the above map for $\nu$ sufficiently large. This 
implies that the above map is surjective and hence bijective.

\proofend


\subsection{Correspondence of chain complexes}

\smallskip
Fix $\varepsilon > 0$ small enough so that the conclusion of Theorem \ref{correspondence}
holds. If we identify $\mathcal{M}^c(q,p)$ with $\mathcal{M}_{h_\varepsilon}(q,p) \times \{0\}$
using Theorem \ref{correspondence}, then 
$$
\mathcal{M}_{h_\varepsilon}(q,p) \times [0,\varepsilon]
$$
determines a trivial smooth cobordism between $\mathcal{M}^c(q,p)$ and
$\mathcal{M}_{h_\varepsilon}(q,p) \approx \mathcal{M}_{h_\varepsilon}(q,p) \times 
\{\varepsilon\}$. If we choose orientations for the unstable manifolds of $h_\varepsilon$,
then $\mathcal{M}_{h_\varepsilon}(q,p)$ becomes an oriented zero dimensional manifold and
there is an induced orientation on $\mathcal{M}_{h_\varepsilon}(q,p) \times [0,\varepsilon]$.

\begin{definition}
Let $p,q \in \text{Cr}(h_\varepsilon)$ with $\lambda_q - \lambda_p = 1$,
define an orientation on the zero dimensional manifold $\mathcal{M}^c(q,p)$ by identifying
it with the left hand boundary of $\mathcal{M}_{h_\varepsilon}(q,p) \times [0,\varepsilon]$.
\end{definition}

\noindent
An orientation on $\mathcal{M}^c(q,p)$ assigns a $+1$ or $-1$ to each point in 
$\mathcal{M}^c(q,p)$. This determines an integer $n^c(q,p) = \#\mathcal{M}^c(q,p)\in \mathbb{Z}$.
Moreover, the one dimensional manifold $\mathcal{M}_{h_\varepsilon}(q,p) \times 
[0,\varepsilon]$ consists of finitely many closed intervals where the right hand boundary 
is identified with $\mathcal{M}_{h_\varepsilon}(q,p)$. Thus,
$$
n^c(q,p) = -n_{h_\varepsilon}(q,p).
$$

\begin{definition}
Define the $k^{\text{th}}$ chain group $C_k^c(f)$ to be the free abelian group generated by the
critical points of total index $k$ of the Morse-Smale functions $f_j$ for all $j=1,\ldots ,l$,
and define $n^c(q,p)$ to be the number of flow lines with cascades between 
a critical point $q$ of total index $k$ and a critical point $p$ of total index $k-1$
counted with signs determined by the orientations. Let
$$
C^c_\ast(f) = \bigoplus_{k=0}^m C_k^c(f)
$$
and define a homomorphism $\partial_k^c:C^c_k(f) \rightarrow C^c_{k-1}(f)$ by
$$
\partial^c_k(q) = \sum_{p \in Cr(f_{k-1})} n^c(q,p)p.
$$
\end{definition}

\begin{corollary}[Correspondence of Chain Complexes]\label{chaincorrespondence}
For $\varepsilon > 0$ sufficiently small we have $C^c_k(f) = C_k(h_\varepsilon)$ and
$\partial^c_k = - \partial_k$ for all $k=0,\ldots ,m$, where $\partial_k$ denotes the
Morse-Smale-Witten boundary operator determined by the Morse-Smale function $h_\varepsilon$.
In particular, $(C_\ast^c(f), \partial^c_\ast)$ is a chain complex whose homology is
isomorphic to the singular homology $H_\ast(M;\mathbb{Z}$).
\end{corollary}

\bibliographystyle{plain}
\bibliography{books,papers}

\end{document}